\documentclass[preprint,12pt]{elsarticle}
\usepackage{lineno}
\usepackage[breaklinks]{hyperref}
\usepackage{amssymb}
\usepackage{subcaption}
\usepackage{ntheorem}
\usepackage{amsmath}
\usepackage{multirow,booktabs}
\theorembodyfont{\upshape}
\theoremseparator{.}

\newtheorem{theorem}{Theorem}[section]
\newtheorem{remark}{Remark}[section]
\newtheorem{proposition}{Proposition}[section]
\newtheorem{lemma}{Lemma}[section]
\newtheorem{definition}{Definition}[section]
\newtheorem{corollary}{Corollary}
\newtheorem*{pro}{Proof}
\numberwithin{equation}{section}

\begin{document}
\begin{frontmatter}


\title{Global boundedness and Allee effect for a nonlocal time fractional reaction-diffusion equation \tnoteref{t1}}
\tnotetext[t1]{This work is supported by the State Key Program of National Natural Science of China under Grant No.91324201. This work is also supported by the Fundamental Research Funds for the Central Universities of China under Grant 2018IB017, Equipment Pre-Research Ministry of Education Joint Fund Grant 6141A02033703 and the Natural Science Foundation of Hubei Province of China under Grant 2014CFB865.}


\author[mymainaddress]{Hui Zhan}\ead{2432593867@qq.com}

\author[mymainaddress]{Fei Gao\corref{cor1}}\ead{gaof@whut.edu.cn}

\author[mymainaddress]{Liujie Guo}\ead{2252987468@qq.com}

\cortext[cor1]{Corresponding author.}

\address[mymainaddress]{Department of Mathematics and Center for Mathematical Sciences, Wuhan University of Technology, Wuhan, 430070, China}

\begin{abstract}
The global boundedness and asymptotic behavior are investigated for the solutions of a nonlocal time fractional reaction-diffusion equation (NTFRDE)
 $$ \frac{\partial^{\alpha }u}{\partial t^{\alpha }}=\Delta u+\mu u^{2}(1-kJ*u)-\gamma u, \qquad(x,t)\in\mathbb{R}^{N}\times(0,+\infty)$$
	with $0<\alpha <1,\beta, \mu ,k>0,N\leq 2$ and $u(x,0)=u_{0}(x)$. Under appropriate assumptions on $J$ and the property of time fractional derivative, it is proved that for any nonnegative and bounded initial conditions, the problem has a global bounded classical solution if $k^{*}=0$ for $N=1$ or $k^{*}=(\mu C^{2}_{GN}+1)\eta^{-1}$ for $N=2$, where $C_{GN}$ is the constant in Gagliardo-Nirenberg inequality. With further assumptions on the initial datum, for small $\mu$ values, the solution is shown to converge to $0$ exponentially or locally uniformly as $t \rightarrow \infty$, which is referred as the Allee effect in sense of Caputo derivative. Moreover, under the condition of $J \equiv 1$, it is proved that the nonlinear NTFRDE has a global bounded solution in any dimensional space with the nonlinear diffusion terms $\Delta u^{m}\, (2-\frac{2}{N}< m\leq 3)$.
\end{abstract}


\end{frontmatter}


\section{Introduction}
 In this work we study the nonlocal  time fractional reaction-diffusion equation
\begin{eqnarray}
\frac{\partial^{\alpha }u}{\partial t^{\alpha }}&=&\Delta u+\mu u^{2}(1-kJ*u) -\gamma u,\quad(x,t)\in\mathbb{R}^{N}\times(0,+\infty),\label{1}\\
u(x,0)&=&u_{0}(x),\quad x\in \mathbb{R}^{N}, \label{2}
	\end{eqnarray}
	with $0< \alpha <1,N\leq 2,\mu ,k,\gamma >0$. By \cite{0Global} , Here $J(x)$ is a competition kernel  with
\begin{equation}\label{3}
  0\leq J \in L^{1}(\mathbb{R}^{N}),\int_{\mathbb{R}^{N}}J(x)dx=1,\inf\limits_{B(0,\delta_{0})} J> \eta
  \end{equation}
   for some $\delta _{0}>0,\eta>0$. where $B(0,\delta_{0})=(-\delta_{0},\delta_{0})^{N}$, $$J*u(x,t)=\int_{\mathbb{R}^{N}}J(x-y)u(y,t)dy,$$
and $F(u,J*u)=\mu u^{2}(1-kJ*u)-\gamma u$ is represented as sexual reproduction in the population dynamics system, where $-\gamma u$ is population mortality. $J*u$ describes the resource consumption of individuals in a certain area around point $x$ at that point, where $J(x-y)$ represents the probability density function, which describes the distribution of available resources near their individual locations. This type of nonlocal model is not only used in biological phenomenon models such as emergence and evolution of biological species and the process of speciation in \cite{0Global}, but its $\mu u^{\alpha }(1-kJ*u^{\beta }),\alpha ,\beta \geq 1$ form is used to describe the selection process in the mutation process in \cite{2019Global}. Both of these two documents give very complete sufficient conditions for the proof of the global boundedness of the nonlocal reaction-diffusion equation. From \cite{0Global}, we find it used the comparison principle and Duhamel's formulas and other methods to prove the case of $\alpha=2 ,\beta=1$ in \cite{2019Global}.

 In the past few years, the nonlocal reaction-diffusion equation has been widely used in biological models, medicine and other fields. For example, in \cite{2017Wavefronts}, in order to solve the population dynamic problem with nonlocal resource consumption, the existence of the solution of the nonlocal reaction-diffusion wavefront equation was proposed. In \cite{2014Existence}, the existence of a solution in the form of a stable pulse is proved by the perturbation method. Through \cite{2019Mathematical}, we learn that the nonlocal advection model has received considerable attention. It solves the problem of how the dynamic model of cells and cell matrix adhere to the tissue rearranges themselves. Such equation can not only be used for cell model research, but also can be used to explain uninfected viruses and cancer cells through the PDE-ODE equation \cite{2020Global}. From \cite{2011Population}, two numerical methods are proposed to predict biological populations, local competition and mutation. Compared with the traditional reaction-diffusion equation, the nonlocal reaction-diffusion equation has new mathematical characteristics and richer nonlinear dynamics. Therefore, many researches are devoted to the nonlocal reaction-diffusion equation. In \cite{2015Preface}, in order to describe the distribution of population density in the case of nonlocal consumption of resources, the following nonlocal reaction-diffusion equation is obtained.
$$\frac{\partial u}{\partial t}=D\frac{\partial ^{2}u}{\partial x^{2}}+F(u,J(u)),$$
where $F(u,J(u))=au^{k}(1-J(u))-\sigma u, J(u)=\int_{-\infty }^{+\infty }\phi (x-y)u(y,t)dy.$ Here $k$ is a positive integer, $k=1$ corresponds to asexual and $k=2$ to sexual reproduction. Then $u = 1$ is the stationary solution of this equation. It is stable in the case of local equations, but may lose stability to nonlocal equations. Therefore, in \cite{2015Preface}, it can be obtained that the local equations have stable monotonic waves and unstable pulses, while the nonlocal equations have simple periodic waves and stable or unstable pulses. As discussed above, the stable solution exists in a bistable situation with global resource consumption. From the above discussion, we can find that in order to study the dynamics of the oncolytic virus model, construct a mathematical model of cell migration, and describe the emergence and evolution of biological species, the above problems consider the existence, stability and instability and other properties of solutions using nonlocal reaction-diffusion equations to solve.

In the nonlocal Fisher-KPP equations in recent years, for the existence and stability of the frontier, the classical solution depends to a large extent on the application of the comparison principle in
\cite{0Global}\cite{2019Global}. But in the following the lack of comparison principle in the local Fisher-KPP equation is one of the characteristics by \cite{2009The}.
$$\frac{\partial u}{\partial t}=\frac{\partial ^{2}u}{\partial x^{2}}+u(1-J_{\sigma }*u),\quad (0,\infty )\times \mathbb{R}$$
where $J_{\sigma }$ is a given convolution kernel and
$$J_{\sigma }*u\left ( x \right ) =\int _{\mathbb{R}^{d} }u(x-y)J\left ( y \right )dy.$$
The equation admits a semi-wavefront connecting $0$ to an unknown positive state for all $c \geq c^{*}= 2$ and there is no such kind of wavefront with wave speed $c < 2$. As mentioned in \cite{1936The}, introducing nonlocal intraspecific competition for resources changes the properties of solutions of this equation. Many studies are devoted to the below nonlocal reaction-diffusion equation	
$$ \frac{\partial u}{\partial t}=\Delta u+u(1-u).$$
Some progress has recently been attained in this direction for the so called nonlocal Fisher-KPP equation
$$\frac{\partial u}{\partial t} =\Delta u+\mu u\left ( 1-J*u \right ), $$
for which $u=1$ is a stationary solution. Comparing the above two nonlocal Fisher-KPP equations, we can find that the properties brought about by different nonlocal terms are also different. The former mainly requires proof of the existence of solutions to nonlocal equations of traveling wave solutions.  This similar proof method is mentioned in \cite{Volpert2015Pulses}\cite{2014Elliptic}.  The latter is stable in the case of local equations, but may lose stability to nonlocal equations. This phenomenon is not only observed in the study of traveling wave solutions, but also has similar results for the following nonlinear nonlocal bistable reaction-diffusion equation
$$\frac{\partial u}{\partial t}=\frac{\partial ^{2}u}{\partial t^{2}}+u^{2}(1-J_{\sigma }*u)-du,\qquad (t,x)\in  (0,\infty )\times \mathbb{R}$$
are investigate in previous work \cite{2017Wavefronts}, where $-du$ is the mortality term and $d$ is the death rate. The reaction term $ u^{2}(1-J_{\sigma }*u)-du$ consists of the reproduction which is proportional to the square of the density, the available resources and the mortality.

Cancer is a disease involving abnormal changes in cell populations, which can cause significant morbidity and mortality. One of the earliest and most prosperous methods in theoretical cancer research is to use available event statistics and create models to explain the observations, see \cite{2018M}. In order to make outstanding contributions to the development of anti-cancer therapy, a large number of mathematical tools need to be introduced to perform different types of models in cancer biology. In the process, time fractional reaction-diffusion equation is used as a modeling tool for many scientific and engineering applications. On the other hand, fractional differential equations have attracted great interest from most people in recent years. Favored by the research community due to its wide range of applications; for example, see previous research \cite{Zhou2017Attractivity}\cite{2016Weak} and \cite{2000Applications} can get some theories about its application. It can be seen that the time fractional reaction-diffusion equation is usually used to establish this type of cell model.

In addition to the practical application of the time fractional reaction-diffusion equation mentioned above, it has many applications. For example, through the deduced maximum regularity of the fractional-order linear reaction-diffusion equation, the existence of the global time solution of a class of time fractional diffusion system with time fractional derivative is studied \cite{2017Global}. In \cite{2015Time}, we propose a definition of the Caputo fractional derivative on a finite interval in the fractional sobolev spaces and investigate it from the operator theoretic viewpoint. In particular, some important equivalences of the norms related to the fractional integration and differentiation operators in the fractional sobolev spaces are given. It deals with time fractional reaction-diffusion equatioan
$$\partial _{t}^{\alpha }u(x,t)=-Lu(x,t)+F(x,t),\quad x\in \Omega \subset \mathbb{R}^{n},0<t\leq T$$
	where $L$ is a differential operator of elliptic type and $\partial _{t}^{\alpha }$ denotes the left Caputo fractional derivative that is usually defined by the formula
\begin{equation}\label{4}
  \partial _{t}^{\alpha }u(x,t)=\frac{1}{\Gamma (1-\alpha )}\int_{0}^{t}\frac{\partial u}{\partial s}(x,s)(t-s)^{-\alpha }ds,\quad 0<\alpha<1\quad
\end{equation}
in the formula \eqref{4}, the left Caputo fractional derivative $\partial _{t}^{\alpha }u$ is a derivative of the order $\alpha,0<\alpha<1$ by \cite{2018M}. Moreover, from the analysis of \cite{2015On}, the following reaction-diffusion equation with a Caputo fractional derivative in time and with various boundary conditions is considered. $$^{c}\textrm{D}^{\alpha }u=\Delta u-u(1-u), x\in \Omega, t>0$$
with the initial condition $u(x,0)=u_{0}(x), x\in\Omega$. the initial data $u_{0}(x)$ is a given positive and bounded function. Here $d>0$ is the diffusion coefficient and $^{c}\textrm{D}^{\alpha }$ is the time fractional derivative of order $0<\alpha<1$ defined by $$^{c}\textrm{D}^{\alpha }u(x,t)=\frac{1}{\Gamma (1-\alpha )}\int_{0}^{t}\frac{\partial u}{\partial s}(x,s)(t-s)^{-\alpha }ds$$
for a function $u$ differentiable in the time variable
\cite{Kilbas2006}\cite{2014Basic}. In \cite{2015On}, it proved the existence of global bounded solutions as well as blowing-up solutions according to the condition imposed on the initial data. Based on the above time fractional reaction-diffusion equation, in order to better prove it, we will define the same the time fractional derivative. From the practical application of the above mentioned time fractional reaction-diffusion equation, we not only get the definition of time fractional derivative, but also understand the existence and uniqueness of the solution of such type of equation.

In recent years, the use of time fractional reaction-diffusion equations in mathematical models has become more and more popular. It solves the numerical resolution of the time fractional reaction-diffusion equation by constructing and analyzing stable high-order schemes by \cite{2007Finite}. A numerical algorithm was proposed in \cite{2016Error}, which uses a high-order finite difference method to solve the time fractional reaction-diffusion equation. In \cite{2018A}, in order to propose different solutions to different time fractional reaction-diffusion equations, it proposed explicit and implicit difference methods for time fractional reaction-diffusion equations. In addition to the research on the existence of the solution to the time fractional reaction-diffusion equation, the nature of the solution to the equation is also involved. For example, the maximum regularity result is used to prove the global boundedness and asymptoticity of the time fractional reaction-diffusion equations \cite{2017Global}. It can be seen that the application of the time fractional reaction-diffusion equation in the mathematical model is of great significance.

From our previous extensive research on the nonlocal reaction-diffusion equation, we found that nonlocality has a profound impact on the reaction-diffusion equation. Therefore, we consider the nonlocal nature of the time fractional reaction-diffusion equation. In \cite{2016Inverse} \cite{2018M}, the nonlocal concept under the time fractional reaction-diffusion equation is mentioned. The former solves the boundary conditions of the inverse problem of the NTFRDE, and the latter uses Faedo-Galerkin to establish the existence and uniqueness of the weak solution of the proposed model to predict tumor invasion and growth. It is worth mentioning that some authors dedicated to the study of time fractional partial differential equations and their applications will use energy inequalities to establish the existence of solutions to boundary value problems, or use energy inequalities to solve time fractional reaction-diffusion wave equations in \cite{Alikhanov2010A}. In \cite{2009Boundary}\cite{2014Asymptotic}, we often try to analyze the initial-boundary value problem of generalized distribution order time fraction reaction-diffusion equation to obtain open and bounded multidimensional domains. Using the maximum principle of the time fractional reaction-diffusion equation of the appropriate generalized distribution order, the uniqueness and existence results of the initial-boundary value problem of this type of equation can be derived \cite{2010Some}\cite{2009Maximum}. The above papers not only put forward the concept of nonlocal, but also put forward the method of solving the time fractional reaction-diffusion equation accordingly. The conditions and proof methods given are also different, such as Laplace transform, maximum value principle.

Next, we will look at the specific expression form and processing method of the time fractional reaction-diffusion equation in practical application. In \cite{2005Time}, by using Mellin transform, Laplace transform, and the properties of fox function to derive the complete solution of the following time fractional reaction-diffusion equation in half space
\begin{equation}
\begin{cases}
 \nonumber &\frac{\partial ^{\alpha }u(x,,t)}{\partial t^{\alpha }}=b^{2}\frac{\partial ^{2}u(x,t)}{\partial x^{2}}-c^{2}u(x,t)\\
\nonumber&u(x,0)=f(x),\qquad  t>0,x\geq 0
\end{cases}
\end{equation}
which is $b^{2}>0,c^{2}>0,0<\alpha \leq 1$, $\frac{\partial ^{\alpha }u(x,t)}{\partial t^{\alpha }}$ is Caputo fractional derivative.The definition of the Caputo fractional derivative in this equation may have changed nowadays, but the method used the use of Laplace transform to prove the existence of the solution is also instructive in the paper. Through \cite{2013Uniqueness}, we understand that its basic model adopts the form of time fractional reaction-diffusion equation, and use the maximum principle of fractional differential equations to propose a numerical algorithm to prove its uniqueness,and its initial boundary value problem with Neumann boundary conditions for the time fractional reaction-diffusion equation in the form.
\begin{equation}
\begin{cases}
 \nonumber\partial _{t}^{\alpha }u(x,t)=\Delta u(x,t)+f(u(x,t)),&0<\alpha <1,x\in \Omega,0<t<T\\
\nonumber \partial _{v}u(x,t)=g(x,t),&x\in \partial \Omega,0<t<T\\
\nonumber u(x,0)=u_{0},&x\in \Omega
\end{cases}
\end{equation}
which is $\Omega \subset \mathbb{R}^{n}$ is a bounded domain with a smooth boundary $\partial \Omega ,v(x)=(v_{1}(x),\cdots ,v_{n}(x))$ denotes the unit outward normal vector to $\partial \Omega$ at $x,\partial _{v}u=\bigtriangledown u\cdot v$. It proves the uniqueness of the solution to the inverse problem of the fractional differential equations, and a numerical algorithm is given for this problem. The definition of the Caputo fractional derivative defined in it is consistent with the paper.

In the paper, we put forward the concept of NTFRDE, in which sufficient conditions are given for the global boundedness of its solution, asymptoticity and global boundedness under nonlinear conditions. The key step in the process of proof is to solve the fractional differential equations. Fractional differential equations which are differential equations of arbitrary order have attracted considerable attention recently \cite{1993An}\cite{1974The}\cite{1991Relaxation}\cite{2002Analytical}\cite{2009Existence}\cite{2010Analytic}.
However, we find that the existence of Laplace transform is taken for granted in some papers to solve fractional differential equations (see \cite{2007Global}\cite{2008Stability}).
In the process of verifying, In order to better scale the solutions of fractional differential equations, we use fractional Duhamel's formula mentioned in \cite{2018Duhamel}. In addition to solving problems using fractional differential equations, we also emphasize and use the forms of their inequalities. In \cite{2011Laplace}, Giving a sufficient condition to guarantee the rationality of solving constant coefficient fractional differential equations by the Laplace transform method. Next, we study the global boundedness and asymptotic for the solution of NTFRDE of \eqref{1}-\eqref{4}. The main results of the paper are the following.
\begin{theorem}\label{theorem1}
Suppose $u:[0,\infty )\times \mathbb{R}^{n}\rightarrow \mathbb{R}$ and $n\geq 2$,the Caputo fractional derivative with respect to time $t$ of $u$ is defined by \eqref{32}. Then there is
\begin{equation}\label{5}
   u^{n-1}(_{0}^{C}\textrm{D}_{t}^{\alpha }u) \geq \frac{1}{n}(_{0}^{C}\textrm{D}_{t}^{\alpha }u^{n}).
\end{equation}
\end{theorem}

\begin{theorem}\label{theorem2}
Suppose $u:[0,\infty )\times\mathbb{R}^{n}\rightarrow \mathbb{R}$,the left Caputo fractional derivative with respect to time $t$ of $u$ is defined by \eqref{4}. Then there is
\begin{equation}\label{6}
  \int _{B(x,\delta )}\partial _{t}^{\alpha }udy=\partial _{t}^{\alpha }\int _{B(x,\delta )}udy.
\end{equation}
	\end{theorem}
\begin{corollary}\label{corollary2}
  Suppose $u:[0,\infty )\times\mathbb{R}^{n}\rightarrow \mathbb{R}$,the Caputo fractional derivative with respect to time $t$ of $u$ is defined by \eqref{32}. Then there is
\begin{equation}\label{77}
  \int _{B(x,\delta )}(_{0}^{C}\textrm{D}_{t}^{\alpha }u)dy=_{0}^{C}\textrm{D}_{t}^{\alpha }\int _{B(x,\delta )}udy.
\end{equation}
\end{corollary}

\begin{theorem}\label{theorem4}
  Suppose \eqref{3}-\eqref{4} holds and $0\leq u_{0}\in L^{\infty }(\mathbb{R}^{N})$ and $T$ is constant. Denote $k^{*}=0$ for $N=1$ and $k^{*}=(\mu C_{GN}^{2}+1)\eta ^{-1}$ for $N=2$, $C_{GN}$ is the constant appears in Gagliardo-Nirenberg inequality in Lemma \ref{lemma2.1}, then for any $k>k^{*}$, the nonnegative solution of \eqref{1}-\eqref{2} exists and is globally bounded in time, that is, there exist
  \begin{equation}
 \nonumber   K=\begin{cases}
 K(\Arrowvert u_{0}\Arrowvert_{L^{\infty}(\mathbb{R}^{N})},\mu,\eta,k,C_{GN},T^{\alpha}),&  N=1,\\
 K(\Arrowvert u_{0}\Arrowvert_{L^{\infty}(\mathbb{R}^{N})},\mu,C_{GN},T^{\alpha }),&  N=2,
      \end{cases}
  \end{equation}
  such that
\begin{equation}\label{7}
     0\leq u(x,t)\leq K,\qquad\forall (x,t)\in \mathbb{R} ^{N}\times [0,\infty).
  \end{equation}
\end{theorem}

\begin{theorem}\label{theorem5}
  Denote $u(x,t)$ the globally bounded solution of \eqref{1}-\eqref{2}

  (i)For any $\gamma>0$, there exist $\mu ^{*}>0$ and $m^{*}>0$ such that for\\ $\mu \in (0,\mu ^{*})$ and $\left \| u_{0} \right \|_{L^{\infty }(\mathbb{R}^{N})}< m^{*}$, we have $$\left \| u(x,t) \right \|_{L^{\infty }(\mathbb{R}^{N}\times [0,+\infty ))}< \frac{\tau }{\mu },$$ and thus $$\left \| u(x,t) \right \|_{L^{\infty }(\mathbb{R}^{N})}\leq \left \| u_{0} \right \|_{L^{\infty }(\mathbb{R}^{N})}e^{\left ( -\sigma  \right )^{\frac{1}{\alpha }}t},$$ for all $t>0,0<\alpha<1$ with $\sigma :=\gamma -\mu \left \| u(x,t) \right \|_{L^{\infty }(\mathbb{R}^{N}\times [0,+\infty ))}>0.$

  (ii)If $0< \gamma < \frac{\mu }{4k} $, there exist $\mu ^{**}>0$ and $m^{**}>0$ such that for $\mu \in (0,\mu ^{**})$ and $\left \| u_{0} \right \|_{L^{\infty }(\mathbb{R}^{N})}< m^{**}$, we have $$\left \| u(x,t) \right \|_{L^{\infty }(\mathbb{R}^{N}\times [0,+\infty ))}< a,$$ and thus $$\lim_{t\rightarrow \infty }u(x,t)=0$$ locally uniformly in $\mathbb{R}^{N}$.
\end{theorem}

\begin{theorem}\label{theorem6}
 \begin{eqnarray}
  \frac{\partial^{\alpha } u}{\partial t^{\alpha }}&=&\bigtriangledown ((u+1)^{m-1}\bigtriangledown u)+u^{2}(1-\int _{\mathbb{R}^{N}}udx)-u(x,t)  \label{1.1.4}\\
u(x,0)&=&u_{0}(x) \qquad x \in \mathbb{R}^{N}\label{1.1.5}
 \end{eqnarray}
with $ (x,t)\in \mathbb{R}^{N}\times(0,+\infty)$. If $2-\frac{2}{N}< m\leq 3,0< \alpha <1$, then the solution to any initial value $0\leq u_{0}\in L^{\infty }(\mathbb{R}^{N})$ exists and \eqref{1.1.4}-\eqref{1.1.5} is globally bounded, exist $M>0$ such that$$0\leq u(x,t)\leq M, \qquad (x,t)\in \mathbb{R}^{N}\times [0,+\infty).$$
	\end{theorem}

\begin{remark}
  The new contribution of this paper lies in the generalization of the framework in studying problems on bounded domain to that of NTFRDE. For the proof of global boundedness, there are three key steps in the paper. Firstly, solved the general solution problem of fractional differential equations by using Laplace transform. Secondly, based on Duhamel's formula in \cite{2018Duhamel}\cite{2013On}\cite{2019Well}, by using some special ways, it is proved that the fractional Duhamel's formula and give some properties of soltion operator. Finally, this paper uses the fractional differential inequality and $L^{r}$ estimation, which can be used to solve the global boundedness of nonliner NTFRDE.
\end{remark}

\section{Global bounded of solutions for a NTFRDE}

\begin{lemma}\textsuperscript{\cite{Alikhanov2010A}}\label{lemma23}
The left Caputo derivative with respect to time $t$ of $v$ is defined by \eqref{4} and absolutely continuous on $[0,T]$, one has the inequality
  \begin{equation}\label{888}
    v(t)\partial _{t}^{\alpha }v(t)\geq \frac{1}{2}\partial _{t}^{\alpha }v^{2}(t),\quad 0<\alpha <1.
  \end{equation}
\end{lemma}

\noindent\textbf{Proof of Theorem \ref{theorem1}}:
we know this relation \eqref{5} is true for $n=2$, let us rewrite inequality \eqref{5} in the form
$$
\begin{aligned}
 &v(t)(_{0}^{C}\textrm{D}_{t}^{\alpha }v(t))-\frac{1}{2}(_{0}^{C}\textrm{D}_{t}^{\alpha }v^{2}(t))\\
 &=\frac{1}{\Gamma (1-\alpha )}v(t)\int_{0}^{t}\frac{{v}'(\tau )d\tau }{(t-\alpha )^{\alpha }}-\frac{1}{2\Gamma (1-\alpha )}\int_{0}^{t}\frac{2v(\tau ){v}'(\tau )d\tau }{(t-\tau )^{\alpha }}\\
&=\frac{1}{\Gamma (1-\alpha )}\int_{0}^{t}\frac{{v}'(\tau )(v(t)-v(\tau ))d\tau }{(t-\tau )^{\alpha }}\\
&=\frac{1}{\Gamma (1-\alpha )}\int_{0}^{t}\frac{{v}'(\tau )d\tau }{(t-\tau )^{\alpha }}\int_{\tau }^{t}{v}'(\eta )d\eta \\
&=\frac{1}{\Gamma (1-\alpha )}\int_{0}^{t}{v}'(\eta )d\eta \int_{0}^{\eta }\frac{{v}'(\tau )d\tau }{(t-\tau )^{\alpha }}\equiv I_{1}\geq 0.
\end{aligned}
$$
Therefore, in order to prove the inequality, it suffices to show that the integral $I_{1}$ is nonnegative. The integral $I_{1}$ takes nonnegative value, since
$$
\begin{aligned}
  &\frac{1}{\Gamma (1-\alpha )}\int_{0}^{t}{v}'(\eta )d\eta \int_{0}^{\eta }\frac{{v}'(\tau )d\tau }{(t-\tau )^{\alpha }}\\
  &=\frac{1}{\Gamma (1-\alpha )}\int_{0}^{t}(t-\eta )^{\alpha }\frac{{v}'(\eta )d\eta }{(t-\eta )^{\alpha }}\int_{0}^{\eta }\frac{{v}'(\tau )d\tau }{(t-\tau )^{\alpha }}\\
&=\frac{1}{2\Gamma (1-\alpha )}\int_{0}^{t}(t-\eta )^{\alpha }\frac{d }{d\eta }\left ( \int_{0}^{\eta } \frac{{v}'(\tau )d\tau }{(t-\tau )^{\alpha }}\right )^{2}d\eta \\
&=\frac{\alpha }{2\Gamma (1-\alpha )}\int_{0}^{t}(t-\eta )^{\alpha -1}\left ( \int_{0}^{\eta } \frac{{v}'(\tau )d\tau }{(t-\tau )^{\alpha }}\right )^{2}d\eta \geq 0,
\end{aligned}
$$
assume that it is true  for some $n=k-1$. If $n=k$, and Lagrange's mean value theorem, then the induction hypothesis implies
\begin{equation}
\begin{aligned}
\nonumber &u^{k-1}(t)(_{0}^{C}\textrm{D}_{t}^{\alpha }u(t))-\frac{1}{k}(_{0}^{C}\textrm{D}_{t}^{\alpha }u^{k}(t))\\
&=u^{k-1}(t)\int_{0}^{t}\frac{1}{\Gamma (1-\alpha )}\frac{{u}'(\tau )d\tau }{(t-\tau )^{\alpha }}-\frac{1}{k}\int_{0}^{t}\frac{1}{\Gamma (1-\alpha )}\frac{ku^{k-1}(\tau ){u}'(\tau )d\tau }{(t-\tau )^{\alpha }}\\
\nonumber&=\frac{1}{\Gamma (1-\alpha )}\int_{0}^{t}\frac{{u}'(\tau )}{(t-\tau )^{\alpha }}(u^{k-1}(t)-u^{k-1}(\tau ))d\tau \\
\nonumber&=\frac{1}{\Gamma (1-\alpha )}\int_{0}^{t}\frac{{u}'(\tau )}{(t-\tau )^{\alpha }}{(u^{k-1}(\xi ))}'(t-\tau )d\tau \\
\nonumber&=\frac{1}{(k-1)\Gamma (1-\alpha )}\int_{0}^{t}(t-\tau )^{1-\alpha }({u(\tau )}')^{2}[u^{k-2}(t)-u^{k-2}(\tau )]d\tau\\
\nonumber&\equiv I_{2}.
\end{aligned}
\end{equation}
Established by $n=k-1$, we can get
$$
\begin{aligned}
  &u^{k-2}(t)(_{0}^{C}\textrm{D}_{t}^{\alpha }u(t))-\frac{1}{k-1}(_{0}^{C}\textrm{D}_{t}^{\alpha }u^{k-1}(t))\\
 & =\frac{1}{\Gamma (1-\alpha )}\int_{0}^{t}\frac{{u}'(\tau )}{(t-\tau )^{\alpha }}[u^{k-2}(t)-u^{k-2}(\tau )]d\tau \geq 0.
\end{aligned}
$$
Thus, $I_{2}=u^{k-1}(t)(_{0}^{C}\textrm{D}_{t}^{\alpha }u(t))-\frac{1}{k}(_{0}^{C}\textrm{D}_{t}^{\alpha }u^{k}(t))\geq 0$. The relation \eqref{5} is proved.
\begin{remark}
  The theorem played a very important role in the following proofs. So we prove this theorem by mathematical induction. The proof method for the special form of $n=2$ is also mentioned in \cite{Alikhanov2010A}, and the proof of the form of $n=k$ is proved by the Lagrangian median theorem.
\end{remark}

\noindent\textbf{Proof of Theorem \ref{theorem2}}: From $\partial _{t}^{\alpha }u$ by \eqref{4}, we obtain
$$\partial _{t}^{\alpha }u(x,t)=\frac{1}{\Gamma (1-\alpha )}\int_{0}^{t}\frac{\partial u}{\partial s}(x,s)(t-s)^{-\alpha }ds,\quad 0<t\leq T,\quad 0<\alpha<1,$$
and
\begin{eqnarray}
\nonumber\partial _{t}^{\alpha }\int _{B(x,\delta )}udy&=&\frac{1}{\Gamma (1-\alpha )}\int_{0}^{t}\frac{\partial }{\partial s}\int _{B(x,\delta )}u(x,s)dy(t-s)^{-\alpha }ds\\
\nonumber&=&\frac{1}{\Gamma (1-\alpha )}\int_{0}^{t}\int _{B(x,\delta )}\frac{\partial u}{\partial s}(x,s)dy(t-s)^{-\alpha }ds\\
\nonumber&=&\frac{1}{\Gamma (1-\alpha )}\int_{0}^{t}\int _{B(x,\delta )}\frac{\partial u}{\partial s}(x,s)(t-s)^{-\alpha }dyds\\
\nonumber&=&\int _{B(x,\delta )}\frac{1}{\Gamma (1-\alpha )}\int_{0}^{t}\frac{\partial u}{\partial s}(x,s)(t-s)^{-\alpha }dsdy\\
\nonumber&=&\int _{B(x,\delta )}\partial _{t}^{\alpha }udy.
\end{eqnarray}
The proof of the theorem is completed, and the proof of the Corollary \ref{corollary2} is the same as the above method.

\begin{definition}\textsuperscript{\cite{Kilbas2006}}
Assume that $X$ is a Banach space and let $u:\left [ 0,T \right ]\rightarrow X$. The Caputo fractional derivative operators of $u$ for order $\alpha \in \mathbb{C},(Re(\alpha )>0)$ are defined by
\begin{eqnarray}
_{0}^{C}\textrm{D}_{t}^{\alpha }u(t)&=&\frac{1}{\Gamma (1-\alpha )}\int_{0}^{t}(t-s)^{-\alpha }\frac{d}{ds}u(s)ds,\label{32}\\
\nonumber_{t}^{C}\textrm{D}_{T}^{\alpha }u(t)&=&\frac{-1}{\Gamma (1-\alpha )}\int_{t}^{T}(s-T)^{-\alpha }\frac{d}{ds}u(s)ds,
\end{eqnarray}
 where $\Gamma (1-\alpha )$ is the Gamma function. The above integrals are called the left-sided and the right-sided the Caputo fractional derivatives
\end{definition}


\begin{definition}\textsuperscript{\cite{Kilbas2006}}
  The Mittag-Leffler function in two parameters is defined as
  $$E_{\alpha ,\beta }(z)=\sum_{k=0}^{\infty }\frac{z^{k}}{\Gamma (\alpha k+\beta )},\qquad z\in \mathbb{C}$$
  where $\alpha >0,\beta >0,\mathbb{C}$ denote the complex plane.
  \end{definition}

\begin{lemma}\textsuperscript{\cite{1948The}}\label{lemma288}
   If $0<\alpha<1,\eta>0$, then there is $0\leq E_{\alpha,\alpha }(-\eta )\leq \frac{1}{\Gamma (\alpha )}$. In addition, for $\eta>0$, $E_{\alpha ,\alpha }(-\eta )$ is a monotonically decreasing function.
\end{lemma}

\begin{lemma}\textsuperscript{\cite{1948The}}\label{lemma299}
  If $0<\alpha<1,t>0$, then there is $0<E_{\alpha ,1}(-t)<1$. In addition, $E_{\alpha,1}(-t)$ is completely monotonous.
\end{lemma}

  \begin{lemma}\textsuperscript{\cite{Kilbas2006}}
    Assume $0<\alpha<2$, for any $\beta \in \mathbb{R}$, there is a constant $\mu$ such that $\frac{\pi \alpha }{2}<\mu <min\left \{ \pi ,\pi \alpha \right \}$, then there is a constant $c=c(\alpha ,\beta ,\mu )>0$ , such that$$\left | E_{\alpha ,\beta }(z)\right |\leq \frac{c}{1+\left | z \right |}, \qquad \mu \leq \left | arg(z) \right |\leq \pi.$$
\end{lemma}

\begin{lemma}\textsuperscript{\cite{1981Geometric}}\label{lem1}
  Suppose $b\geq 0,\beta >0$ and $a(t)$ is a nonnegative function locally integrable on $0\leq t<T$ (some $T\leq +\infty$). And suppose $u(t)$ is nonnegative and locally integrable on $0\leq t<T$ with $$u(t)\leq a(t)+b\int_{0}^{t}(t-s)^{\alpha -1}u(s)ds$$ on this interval. Then$$u(t)\leq a(t)+\theta \int_{0}^{t}E_{\beta }^{'}(\theta (t-s)a(s))ds,\qquad 0\leq t<T.$$ where $$\theta =(b\Gamma (\beta ))^{\frac{1}{\beta }},\qquad  E_{\beta }(z)=\sum_{n=0}^{\infty }z^{n\beta }/\Gamma (n\beta +1),\qquad E_{\beta }^{'}=\frac{d}{dz}E_{\beta }(z)$$
 $$E_{\beta }^{'}(z)\simeq z^{\beta -1}/\Gamma (\beta ) \quad as\quad z\rightarrow 0+,\qquad E_{\beta }^{'}(z)\simeq \frac{1}{\beta }e^{z} \quad as\quad z\rightarrow +\infty$$
  and $E_{\beta }(z)\simeq \frac{1}{\beta }e^{z}$ as $z\rightarrow +\infty$. If $a(t)\equiv a$ is constant, then $$u(t)\leq aE_{\beta }(\theta t).$$
\end{lemma}

\begin{lemma}\textsuperscript{\cite{2011Laplace}}
Let $C$ be complex plane, for any $\alpha>0,\beta>0$ and\\ $A\in C^{n\times n}$,$$L\left \{ t^{\beta -1} E_{\alpha ,\beta }(At^{\alpha })\right \}= s^{\alpha -\beta }(s^{\alpha }-A)^{-1}$$ hold for $\Re s>\left \| A \right \|^{\frac{1}{\alpha }}$, where $\Re s$ represents the real part of the complex number $s$.
\end{lemma}
\begin{remark}
Before proving the global boundedness of equation \eqref{1}-\eqref{2}, we use Laplace's transform to solve the fractional differential equations
\begin{equation}\label{211}
\begin{aligned}
  _{0}^{C}\textrm{D}_{t}^{\alpha }w(t)&=Aw(t)+C_{3},\qquad  0<\alpha <1,t\geq 0\\
 w(0)&=\eta,
 \end{aligned}
\end{equation}
where $_{0}^{C}\textrm{D}_{t}^{\alpha }$ is the Caputo fractional derivative operator \eqref{32}, $A,C_{3}$ is all constants.
\end{remark}
\begin{lemma}\label{lemma2.4}
  Assume \eqref{211} has a unique continuous solution $w(t),t \in [0,\infty)$, and $A,C_{3}$ is constants, then $\exists T\in [0,\infty )$ such that  $w(t)$ is exponentially bounded in $t \geq T$, which is
  \begin{equation}
    \left \| w\left ( t \right )  \right\| \leq (\left \| \eta \right \|e^{-\sigma T}+\frac{(K+\left | C_{3} \right |)T^{\alpha }e^{-\sigma T}}{\alpha \Gamma (\alpha )})Ce^{[A^{\frac{1}{\alpha }}+\sigma]t}.
  \end{equation}

\end{lemma}

\begin{pro}
  It is easy to see that Equation \eqref{211} is equivalent to the following integral equation

  \begin{equation}\label{2.3}
   w(t)=\eta +\frac{1}{\Gamma (\alpha )}\int_{0}^{t}(t-\tau )^{\alpha -1}[Aw(\tau)+C_{3}]d\tau ,0\leq t<\infty
  \end{equation}
  for $t\geq T$, \eqref{2.3} can be rewritten as
  \begin{equation}
  \begin{aligned}
  \nonumber w(t)&=\eta +\frac{1}{\Gamma (\alpha )}\int_{0}^{T}(t-\tau )^{\alpha -1}[Aw(\tau)+C_{3}]d\tau\\
  \nonumber&\qquad+\frac{1}{\Gamma (\alpha )}\int_{T}^{t}(t-\tau )^{\alpha -1}[Aw(\tau)+C_{3}]d\tau.
   \end{aligned}
  \end{equation}
 In view of assumption of Lemma \ref{lemma2.4}, the solution $w(t)(w(0)=\eta)$ is unique and continuous on $[0,\infty)$. Then $Aw(\tau)+C_{3}$ is bounded on $[0,T]$, there exists a constant $K>0$ such that $ \left \| Aw(t)+C_{3}\right \|\leq K$. we have
 \begin{equation}
 \begin{aligned}
 \nonumber\left \| w(t) \right \| &\leq\left \| \eta \right \| +\frac{K}{\Gamma (\alpha )}\int_{0}^{T}(t-\tau )^{\alpha -1}d\tau +\frac{ A }{\Gamma (\alpha )}\int_{T}^{t}(t-\tau )^{\alpha -1}\left \| w(\tau) \right \|d\tau\\
 &\qquad+\frac{\left | C_{3} \right |}{\Gamma (\alpha )}\int_{T}^{t}(t-\tau )^{\alpha -1}d\tau\\
\nonumber&\leq \left \| \eta \right \| +\frac{K}{\Gamma (\alpha )}\int_{0}^{T}(t-\tau )^{\alpha -1}d\tau +\frac{A}{\Gamma (\alpha )}\int_{T}^{t}(t-\tau )^{\alpha -1}\left \| w(\tau) \right \|d\tau\\
&\qquad+\frac{(t-T)^{\alpha }\left | C_{3} \right |}{\alpha \Gamma (\alpha )}.
 \end{aligned}
 \end{equation}
 Multiply this inequality by $e^{-\sigma t}$ and note that $e^{-\sigma t}\leq e^{-\sigma T},e^{-\sigma t}\leq e^{-\sigma \tau }$, to obtain
 $$
 \begin{aligned}
 \left \| w(t) \right \|e^{-\sigma t} &\leq \left \| \eta \right \|e^{-\sigma t} +\frac{Ke^{-\sigma t}}{\Gamma (\alpha )}\int_{0}^{T}(t-\tau )^{\alpha -1}d\tau \\
 &\qquad+\frac{ Ae^{-\sigma t}}{\Gamma (\alpha )}\int_{T}^{t}(t-\tau )^{\alpha -1}\left \| w(\tau) \right \|d\tau+\frac{(t-T)^{\alpha }\left | C_{3} \right |}{\alpha \Gamma (\alpha )}e^{-\sigma t}\\
 &\leq \left \| \eta \right \|e^{-\sigma T} +\frac{Ke^{-\sigma t}}{\alpha \Gamma (\alpha )}(t^{\alpha }-(t-T)^{\alpha })\\
 &\qquad+\frac{ Ae^{-\sigma t}}{\Gamma (\alpha )}\int_{T}^{t}(t-\tau )^{\alpha -1}\left \| w(\tau) \right \|d\tau+\frac{(t-T)^{\alpha }\left | C_{3} \right |}{\alpha \Gamma (\alpha )}e^{-\sigma t}\\
 &\leq \left \| \eta \right \|e^{-\sigma T} +\frac{(K+\left | C_{3} \right |)T^{\alpha }e^{-\sigma T}}{\alpha \Gamma (\alpha )}\\
 &\qquad+\frac{A}{\Gamma (\alpha )}\int_{0}^{t}(t-\tau )^{\alpha -1}\left \| w(\tau) \right \|e^{-\sigma \tau}d\tau , t\geq T.
 \end{aligned}
 $$
 Denote
 \begin{equation}\label{2.4}
 a=\left \| \eta \right \|e^{-\sigma T}+\frac{(K+\left | C_{3} \right |)T^{\alpha }e^{-\sigma T}}{\alpha \Gamma (\alpha )},\quad b=\frac{A}{\Gamma \left ( \alpha  \right )},\quad u(t)=\left \| w(t) \right \|e^{-\sigma t}.
 \end{equation}
We get $$ u(t)\leq a+b\int_{0}^{t}(t-\tau )^{\alpha -1}u(\tau )d\tau ,\quad t\geq T.$$
By Lemma \ref{lem1} $$u(t)\leq aE_{\alpha }(\theta t)=a\sum_{n=0}^{\infty }\frac{(b\Gamma \left ( \alpha  \right ))^{n}t^{n\alpha }}{\Gamma (n\alpha +1)},\quad t\geq T$$
recall that the Mittag-Leffler function is defined as$$E_{\alpha }(t)=\sum_{n=0}^{\infty }\frac{t^{n}}{\Gamma (n\alpha +1)},\quad \alpha >0.$$
Then
\begin{equation}\label{2.5}
u(t)\leq aE_{\alpha }(b\Gamma \left ( \alpha  \right )t^{\alpha}),\quad t\geq T.
\end{equation}
It is known that (see \cite{2001Fractional},p.21) the Mittag-Leffler type function $E_{\alpha }(wt^{\alpha })$ satisfies the following inequality
\begin{equation}\label{2.6}
E_{\alpha }(wt^{\alpha })\leq Ce^{w^{\frac{1}{\alpha }t}},\quad t\geq 0,w\geqslant 0,0<\alpha <2
\end{equation}
where $C$ is a positive constant.
With \eqref{2.5} and \eqref{2.6}, we have $$u(t)\leq aCe^{(b\Gamma (\alpha ))^{\frac{1}{\alpha }}t},\quad t\geq T.$$Then, obtain the following form by \eqref{2.4}
\begin{equation}
\begin{aligned}
\nonumber&\left \| w\left ( t \right )  \right\|=u(t)e^{\sigma t} \leq aCe^{[(b\Gamma (\alpha ))^{\frac{1}{\alpha }}+\sigma ]t}\\
&\leq (\left \| \eta \right \|e^{-\sigma T}+\frac{(K+\left | C_{3} \right |)T^{\alpha }e^{-\sigma T}}{\alpha \Gamma (\alpha )})Ce^{[A^{\frac{1}{\alpha }}+\sigma]t}<\infty,\quad t\geq T.
\end{aligned}
\end{equation}
So $w(t)$ is bounded on $t \in [0,\infty)$.
\end{pro}

\begin{remark}
If the constant term $C_{3}$ of equation \eqref{211} is replaced by $f(t)$.\\ From \cite{2011Laplace}, we know the solution of the following fractional differential equation by Lemma \ref{rem1}
\begin{eqnarray}\label{255}
   _{0}^{C}\textrm{D}_{t}^{\alpha }w(t)&=&Aw(t)+f(t),\qquad t\geq 0 \\
  \nonumber w(0)&=&\eta.
\end{eqnarray}
\end{remark}

\begin{lemma}\textsuperscript{\cite{2011Laplace}}\label{rem1}
  Assume \eqref{255} has a unique continue solution $w(t)$. If $f(t)$ is continuous on $[0,\infty)$ and exponentially bounded, then $w(t)$ and its Caputo fractional derivative $ _{0}^{C}\textrm{D}_{t}^{\alpha }w(t)$ are both exponentially bounded, thus its Laplace transform exists and which is
 \begin{equation}\label{2.7}
 w(t)=E_{\alpha ,1}(t)\eta +\int_{0}^{t}(t-\tau)^{\alpha -1}E_{\alpha ,\alpha }(A(t-\tau)^{\alpha})f(\tau)d\tau.
 \end{equation}
\end{lemma}
\begin{remark}
In the process of proving the global boundedness, we firstly will use fractional Doumal's formula. Therefore, we will prove its related lemma next. We define the operator $A$ on $L^{2}(\Omega )$ by
\begin{equation}
\begin{cases}
  \nonumber D(A)=\left \{ u\in H_{0}^{1}(\Omega ):\Delta u\in L^{2}(\Omega ) \right \}, \\
  \nonumber Au=\Delta u,\qquad for \, u\in D(A).
\end{cases}
\end{equation}
From equations \eqref{1}-\eqref{2}, we consider the equations

\begin{equation}\label{2.8}
\begin{cases}
  \frac{\partial ^{\alpha }}{\partial t^{\alpha }}u(t)=Au(t)+h(t), \\
 u(0)=u_{0},
\end{cases}
\end{equation}
 where $h(t)=\mu u^{2}(1-kJ*u)-\gamma u$ and $A$ is self-adjoint on Hilbert space $X$. Then it follows from spectral theorem that there exists a measure space $(\Sigma ,\mu )$ and a Borel measurable function $a$ and a unitary map $U:L^{2}(\Sigma,\mu )\rightarrow X$ such that$$U^{-1}AU=T_{a},T_{a}\varphi (\xi )=a(\xi )\varphi (\xi ),\quad \xi \in \Sigma. $$
\end{remark}
 \begin{lemma}\label{lem3}
   If $u$ satisfies Equation \eqref{2.8}, then $u$ also satisfies
   $$
   \begin{aligned}
   u(t)&=U(E_{\alpha }(a(\xi )t^{\alpha })U^{-1}u_{0})\\
   &\qquad+\int_{0}^{t}(t-s)^{\alpha -1}U(E_{\alpha ,\alpha }((t-s)^{\alpha }a(\xi ))U^{-1}h(s))ds, \quad t\in [0,T].
    \end{aligned}
   $$
 \end{lemma}

 \begin{pro}
   Let $v(t,\xi )=(U^{-1}u(t))(\xi )$, then we have $$(U^{-1}Au(t))(\xi )=U^{-1}AUv(t,\xi )=T_{a}v(t,\xi )=a(\xi )v(t,\xi ).$$Notice that $$U^{-1}(\frac{\partial ^{\alpha }u(t)}{\partial t^{\alpha }})=\frac{\partial ^{\alpha }}{\partial t^{\alpha }}(U^{-1}u(t))=\frac{\partial ^{\alpha }v(t,\xi )}{\partial t^{\alpha }}.$$We apply the map $U^{-1}$ to both sides of Equation \eqref{2.8}, then we obtain the following equation on $L^{2}(\Sigma)$

   \begin{equation}\label{2.9}
   \begin{cases}
     \frac{\partial ^{\alpha }v(t,\xi )}{\partial t^{\alpha }}=a(\xi )v(t,\xi )+U^{-1}h(t), \\
     v(0)=U^{-1}u_{0}.
   \end{cases}
   \end{equation}
   We know from Remark \ref{rem1} that Equation \eqref{2.9} has a unique solution given by$$v(t,\xi )=E_{\alpha }(a(\xi )t^{\alpha })U^{-1}u_{0}+\int_{0}^{t}(t-s)^{\alpha -1}E_{\alpha ,\alpha }((t-s)^{\alpha }a(\xi ))U^{-1}h(s)ds,t\in [0,T]$$which yields that
   $$
   \begin{aligned}
   u(t)=Uv(t,\xi)&=U(E_{\alpha }(a(\xi )t^{\alpha })U^{-1}u_{0})\\
   &\qquad+\int_{0}^{t}(t-s)^{\alpha -1}U(E_{\alpha ,\alpha }((t-s)^{\alpha }a(\xi ))U^{-1}h(s))ds.
   \end{aligned}
   $$
 \end{pro}
By define two maps $$\mathcal{S} _{\alpha }(t)\phi=U(E_{\alpha }(a(\xi )t^{\alpha })U^{-1}\phi),\, \mathcal{K} _{\alpha }(t)\phi=U(E_{\alpha ,\alpha }((t)^{\alpha }a(\xi ))U^{-1}\phi),\, \phi \in X,$$then through Lemma \ref{lem3}, we can get the following fractional Duhamel's formula.

\begin{lemma}\label{lemma2.6}
  If $u\in C([0,T],X)$ satisfies Equation \eqref{1} -\eqref{2}, then $u$ satisfies the following integral equation:
   $$u(t)=\mathcal{S}_{\alpha }(t)u_{0}+\int_{0}^{t}(t-s)^{\alpha -1}\mathcal{K} _{\alpha }(t-s)h(s)ds.$$
   The solution operators $\mathcal{S}_{\alpha }(t)$ and $\mathcal{K} _{\alpha }(t)$ are defined by the functional calculus of $A$ via the Mittag-Leffler function when evaluated at $A$.
\end{lemma}

 \begin{remark}
   By reference \cite{2018Duhamel}, we are well known that fractional Duhamel's formula plays a very important role in the research for time fractional Schrödinger equation on a Hilbert space. However, Based on proof method, We will use the above fractional Duhamel's formula to decompress the solution for NTFRDE.
 \end{remark}

 \begin{proposition}\textsuperscript{\cite{2018Duhamel}}\label{proposition.1}
   For each fixed $t \geq 0$,$\mathcal{S}_{\alpha }(t)$ and $\mathcal{K} _{\alpha }(t)$ are linear and bounded operators, for any $\phi \in X$,

   \begin{equation}\label{9}
     \left \|\mathcal{S}_{\alpha }(t)\phi  \right \|_{X}\leq C_{4}\left \| \phi \right \|_{X},\qquad  \left \|\mathcal{K}_{\alpha }(t)\phi  \right \|_{X}\leq C_{4}\left \| \phi \right \|_{X},
   \end{equation}
where $C_{4}$ is a constant.
 \end{proposition}

\begin{lemma}\textsuperscript{\cite{1966}}\label{lemma2.1}
  Let $\Omega$ be an open subset of $\mathbb{R}^{N}$, assume that $1\leq p,q \leq \infty$ with $(N-q)p<Nq $ and $r\in(0,p)$. Then there exists constant $C_{GN}>0$ only depending on $q, r$ and $\Omega$ such that for any $u\in W^{1.q}(\Omega)\cap L^{p}(\Omega)$
  $$\int_{\Omega}u^{p}dx\leq C_{GN}(\Arrowvert\bigtriangledown u\Arrowvert^{(\lambda^{*}p)}_{L^{q}(\Omega)}\Arrowvert u\Arrowvert_{L^{p}(\Omega)}^{(1-\lambda^{*})p}+\Arrowvert u\Arrowvert_{L^{p}(\Omega)}^{p})$$
  holds with $$\lambda^{*}=\frac{\frac{N}{r}-\frac{N}{p}}{1-\frac{N}{q}+\frac{N}{r}}\in (0,1).$$
\end{lemma}

\begin{lemma}\textsuperscript{\cite{2019Global}}\label{lemma2.8}
  Fix $x_{0}= (x_{1}^{0},\cdots,x_{N}^{0})\in \mathbb{R}^{N}$, choose $0<\delta\leq \frac{1}{2}\delta_{0}$, and denote $B(x_{0},\delta)=(x_{0}= (x_{1}^{0},\cdots,x_{N}^{0})\in \mathbb{R}^{N}|\quad \lvert x_{i}-x^{0}_{i}\lvert\leq\delta,1\leq i\leq N)$, then for any $y\in B(x,\delta)$, there is
  \begin{equation}\label{10}
  \begin{aligned}
    \int_{B(x,\delta)}\Delta u^{2}dy&=\int_{B(0,\delta)}\Delta u^{2}(y+x,t)dy\\
    &=\Delta \int_{B(0,\delta)}u^{2}(y+x,t)dy=\Delta \int_{B(x,\delta)}u^{2}dy.
    \end{aligned}
  \end{equation}
\end{lemma}

\begin{lemma}\textsuperscript{\cite{2016Partial}}\label{Young}
  $\forall a,b\geq 0$ and $\varepsilon  >0$, for $1<p,q<\infty ,\frac{1}{p}+\frac{1}{q}=1$, then there is $$a\cdot b\leq \varepsilon\frac{a^{p}}{p} +\varepsilon ^{-\frac{q}{p}}\frac{b^{q}}{q}.$$
\end{lemma}

\begin{proposition}\textsuperscript{\cite{1966On}}\label{Proposition2.1}
  Assume the initial data $0\leq u_{0}\in L^{\infty}(\mathbb{R}^{N})$. Then there are a maximal existence time $T_{max}\in(0,\infty]$ and $u\in C([0,T_{max}),L^{\infty}(R^{N}))\bigcap$\\ $C^{2,1}(\mathbb{R}^{N}\times[0,T_{max}))$ such that $u$ is the unique nonnegative classical solution of \eqref{1}-\eqref{2}. Furthermore, if $T_{max}<+\infty$, then$${\lim_{t\to T_{max}}}\Arrowvert u(\cdot,t)\Arrowvert_{L^{\infty}(\mathbb{R}^{N})}=\infty.$$
\end{proposition}

\noindent \textbf{Proof of Theorem \ref{theorem4}}.

Fix $x_{0}= (x_{1}^{0},\cdots,x_{N}^{0})\in \mathbb{R}^{N}$, choose $0<\delta\leq \frac{1}{2}\delta_{0}$, and denote $B(x_{0},\delta)=(x_{0}= (x_{1}^{0},\cdots,x_{N}^{0})\in \mathbb{R}^{N}|\quad \lvert x_{i}-x^{0}_{i}\lvert\leq\delta,1\leq i\leq N)$ with $\lvert B(x_{0},\delta)\lvert=(2\delta)^{N}$. For any $x \in \mathbb{R}^{N}$,multiply \eqref{1} by $2u\varphi_{\varepsilon}$, where $\varphi_{\varepsilon}(\cdot)\in C_{0}^{\infty}(B(x,\delta))$, and $\varphi_{\varepsilon}(\cdot)\to 1$ locally uniformly in $B(x,\delta)$ as $\varepsilon\to 0$. Integrating by parts over $B(x,\delta)$, we obtain

$$
\begin{aligned}
\int_{B(x,\delta)}2u\varphi_{\varepsilon}\frac{\partial^{\alpha } u}{\partial t^{\alpha }}dy&=\int_{B(x,\delta)}2u\varphi_{\varepsilon}\Delta udy+\int_{B(x,\delta)}2u^{3}\mu\varphi_{\varepsilon}(1-kJ*u)dy\\
&\qquad-2\gamma\int _{B(x,\delta )}u^{2}\varphi_{\varepsilon} dy,
\end{aligned}
$$
from $$\int_{B(x,\delta)}2u\varphi_{\varepsilon}\Delta udy=\int_{B(x,\delta)}\Delta u^{2}\varphi_{\varepsilon}dy-2\int_{B(x,\delta)}\arrowvert\bigtriangledown u\arrowvert^{2}\varphi_{\varepsilon}dy,$$
and Theorem \ref{theorem2} and Lemma \ref{lemma23}, then we can get
$$\int_{B(x,\delta)}2u\varphi_{\varepsilon}\frac{\partial^{\alpha } u}{\partial t^{\alpha }}dy\geq\frac{\partial^{\alpha }}{\partial t^{\alpha }}\int_{B(x,\delta)}u^{2}\varphi_{\varepsilon}dy,$$
and
$$
\begin{aligned}
   &\frac{\partial^{\alpha }}{\partial t^{\alpha }}\int_{B(x,\delta)}u^{2}\varphi_{\varepsilon}dy+2\int_{B(x,\delta)}\arrowvert\bigtriangledown u\arrowvert^{2}\varphi_{\varepsilon}dy \\
   &\leq\int_{B(x,\delta)}\Delta u^{2}\varphi_{\varepsilon}dy+2\mu\int_{B(x,\delta)}u^{3}(1-kJ*u)\varphi_{\varepsilon}dy-2\gamma\int_{B(x,\delta)}u^{2}\varphi_{\varepsilon}dy.
\end{aligned}
$$

\noindent Taking $\varepsilon\to 0$, we obtain
\begin{equation}\label{11}
\begin{aligned}
&\frac{\partial^{\alpha}}{\partial t^{\alpha}}\int_{B(x,\delta)}u^{2}dy+2\int_{B(x,\delta)}\arrowvert\bigtriangledown u\arrowvert^{2}dy\\
&\leq\int_{B(x,\delta)}\Delta u^{2}dy+2\mu\int_{B(x,\delta)}u^{3}(1-kJ*u)dy-2\gamma\int_{B(x,\delta)}u^{2}dy\\
&\leq\int_{B(x,\delta)}\Delta u^{2}dy+2\mu\int_{B(x,\delta)}u^{3}dy-2\mu\int_{B(x,\delta)}u^{3}kJ*udy\\
&\qquad-2\gamma\int_{B(x,\delta)}u^{2}dy.
\end{aligned}
\end{equation}
Knowing from \eqref{3} and Lemma \ref{lemma2.8}, we have used the fact that $\forall y\in B(x,\delta),z\in B(y,\delta)$, then $J(z-y)\geq\eta$ and $B(x,\delta)\subset B(y,2\delta)$, and then
$$
\begin{aligned}
&J*u(y,t)=\int_{B(y,\delta)}J(y-z)u(z,t)dz\geq\eta\int_{B(x,\delta)}u(y,t)dy,\\
&\int_{B(x,\delta)}\Delta u^{2}dy=\Delta \int_{B(x,\delta)}u^{2}dy.
\end{aligned}
$$
Therefore
$$
\begin{aligned}
&\frac{\partial^{\alpha }}{\partial t^{\alpha }}\int_{B(x,\delta)}u^{2}dy+2\int_{B(x,\delta)}\arrowvert\bigtriangledown u\arrowvert^{2}dy\\
&\leq\Delta \int_{B(x,\delta)}u^{2}dy+2\mu\int_{B(x,\delta)}u^{3}dy-2\mu\eta k\int_{B(x,\delta)}u^{3}dy\int_{B(x,\delta)}udy\\
&\qquad-2\gamma\int_{B(x,\delta)}u^{2}dy.
\end{aligned}
$$
Now we proceed to estimate the term $\int_{B(x,\delta)}u^{3}dy$. Firstly using Gagliardo-Nirenberg inequality in Lemma \ref{lemma2.1} on $\Omega=B(x,\delta)$ with $p=3,q=r=2,\lambda^{*}=\frac{N}{6}$, there exists constant $C_{GN}>0$, such that
\begin{equation}\label{22}
  \int_{B(x,\delta)}u^{3}dy\leq C_{GN}(N,\delta)(\Vert\bigtriangledown u\Vert_{L^{2}(B(x,\delta))}^{\frac{N}{2}}\Vert u\Vert_{L^{2}(B(x,\delta))}^{3-\frac{N}{2}}+\Vert u\Vert_{L^{2}(B(x,\delta))}^{3}).
\end{equation}
On the one hand, by Lemma \ref{Young}, we can make $a=\Vert\bigtriangledown u\Vert_{L^{2}(B(x,\delta))}^{\frac{N}{2}}$,\\ $b=C_{GN}(N,\delta)\Vert u\Vert_{L^{2}(B(x,\delta))}^{3-\frac{N}{2}}, \varepsilon=\frac{1}{\mu},p=\frac{1}{N},q=\frac{4}{4-N}$, then obtain
\begin{equation}\label{23}
\begin{aligned}
&C_{GN}(N,\delta)\Vert\Delta u\Vert_{L^{2}(B(x,\delta))}^{\frac{N}{2}}\Vert u\Vert_{L^{2}(B(x,\delta))}^{3-\frac{N}{2}}\\
&\leq\frac{1}{\mu}\Vert \Delta u\Vert_{L^{2}(B(x,\delta))}^{2}+\mu^{\frac{N}{4-N}}C^{\frac{4}{4-N}}_{GN}(N,\delta)\Vert u\Vert_{L^{2}(B(x,\delta))}^{\frac{2(6-N)}{4-N}}.
\end{aligned}
\end{equation}
Here by Young's inequality, let $a=C_{GN}(N,\delta ),b=\left \| u \right \|_{L^{2}(B(x,\delta ))}^{3},p=\frac{2(6-N)}{3(4-N)}$, $q=\frac{2(6-N)}{N}$, then
\begin{equation}\label{24}
  C_{GN}(N,\delta )\left \| u \right \|_{L^{2}(B(x,\delta ))}^{3}\leq \left \| u \right \|_{L^{2}(B(x,\delta ))}^{\frac{2(6-N)}{4-N}}+C^{\frac{2(6-N)}{N}}(N,\delta).
\end{equation}
By interpolation inequality, we obtain
\begin{equation}\label{25}
\begin{aligned}
\left \| u \right \|_{L^{2}(B(x,\delta ))}^{\frac{2(6-N)}{4-N}}&\leq (\left \| u \right \|_{L^{1}(B(x,\delta ))}^{\frac{1}{4}}\left \| u \right \|_{L^{3}(B(x,\delta ))}^{\frac{3}{4}})^{\frac{2(6-N)}{4-N}}\\
&=(\int _{B(x,\delta)}u^{3}dy\int _{B(x,\delta )}udy)^{\frac{6-N}{2(4-N)}}.
\end{aligned}
\end{equation}
Combing \eqref{22}-\eqref{25}, we obtain
\begin{equation}
  \begin{aligned}
  &2\mu \int _{B(x,\delta )}u^{3}dy\\
  &\leq 2\int _{B(x,\delta )}\left | \bigtriangledown u \right |^{2}dy+2\mu(1+\mu^{\frac{N}{4-N}}C_{GN}^{\frac{4}{4-N}}(N,\delta ))(\int _{B(x,\delta )}u^{3}dy\int _{B(x,\delta )}udy)^{\frac{6-N}{2(4-N)}}\\
  &\qquad+2\mu C_{GN}^{\frac{2(6-N)}{N}}(N,\delta).\label{26}
  \end{aligned}
\end{equation}
Next we consider the cases $N=1$ and $N=2$ respectively.

\noindent \textbf{Case 1}. $N=1$. From \eqref{26}, by Lemma \ref{Young}, let $\varepsilon =\eta k,p=\frac{6}{5},q=6$, we get
\begin{equation}
  \begin{aligned}
  &2\mu \int _{B(x,\delta )}u^{3}dy\\
&\leq 2\int _{B(x,\delta )}\left | \bigtriangledown u \right |^{2}dy
+2\mu \eta k\int _{B(x,\delta )}u^{3}dy\int _{B(x,\delta )}udy\\
&\qquad+2\mu (\mu^{\frac{1}{3}}C_{GN}^{\frac{4}{3}}(1,\delta )+1)^{6}(\eta k)^{-5}+2\eta C_{GN}^{10}(1,\delta ),\label{27}
  \end{aligned}
\end{equation}
inserting \eqref{27} into \eqref{11}, we have
 \begin{equation}\label{28}
   \begin{aligned}
  &\frac{\partial^{\alpha } }{\partial t^{\alpha }}\int _{B(x,\delta )}u^{2}dy+2\gamma \int _{B(x,\delta )}u^{2}dy\\
&\leq \Delta \int _{B(x,\delta )}u^{2}dy+2\mu [(\mu^{\frac{1}{3}}C_{GN}^{\frac{4}{3}}(1,\delta ))^{6}(\eta k)^{-5}+C_{GN}^{10}(1,\delta )].
  \end{aligned}
 \end{equation}
 Denote $w(t)=\int_{B(x,\delta )}u^{2}dy$, the solution of the following fractional ordinary differential equation
 \begin{equation}\label{29}
   \begin{cases}
    _{0}^{C}\textrm{D}_{t}^{\alpha }w(t)=-2\gamma w(t)+2\mu \left ( (\mu ^{\frac{1}{3}}C_{GN}^{\frac{4}{3}}(1,\delta)+1)^{6}(\eta k)^{-5}+C_{GN}^{10}(1,\delta ) \right )\\
w(0)=2\delta \left \| u_{0} \right \|_{L^{\infty }(R^{N})}^{2}
   \end{cases}
 \end{equation}
 from Lemma \ref{lemma2.4}, let$$A=-2\gamma,C_{3}=2\mu \left ( (\mu ^{\frac{1}{3}}C_{GN}^{\frac{4}{3}}(1,\delta)+1)^{6}(\eta k)^{-5}+C_{GN}^{10}(1,\delta ) \right ),$$
 and $w(0)=\eta =2\delta \left \| u_{0} \right \|_{L^{\infty }(R^{N})}^{2}$, for any $(x,t)\in \mathbb{R}\times [0,T_{max})$, we obtain
 $$
 \begin{aligned}
 &\int _{B(x,\delta )}u^{2}dy=\left \| u \right \|_{L^{2}(B(x,\delta ))}^{2}\leq \left \| w(t) \right \|\\
&\leq w(0)+\frac{(K+\left | C_{3} \right |)T^{\alpha }}{\alpha \Gamma (\alpha )}Ce^{[(-2\gamma )^{\frac{1}{\alpha }}+\sigma]t-\sigma T}\\
&\leq 2\delta \left \| u_{0} \right \|_{L^{\infty }(R^{N})}^{2}+\frac{(K+\left | 2\mu \left ( (\mu ^{\frac{1}{3}}C_{GN}^{\frac{4}{3}}(1,\delta)+1)^{6}(\eta k)^{-5}+C_{GN}^{10}(1,\delta ) \right ) \right |)T^{\alpha }}{\alpha \Gamma (\alpha )}C\\
&:=M_{1}.
 \end{aligned}
 $$
 \noindent \textbf{Case 2}.$N=2$. From \eqref{26}, we have
 \begin{equation}\label{2.10}
   \begin{aligned}
   &2\mu \int _{B(x,\delta )}u^{3}dy\\
   &\leq 2\int _{B(x,\delta )}\left | \bigtriangledown u \right |^{2}dy+2\mu (1+\mu C_{GN}^{2}(2,\delta ))(\int _{B(x,\delta )}u^{3}dy\int _{B(x,\delta )}udy)\\
   &\qquad+2\mu C_{GN}^{4}(2,\delta).
   \end{aligned}
 \end{equation}
 Inserting \eqref{2.10} into \eqref{11}, for $k\geq k^{*}=\frac{\mu C_{GN}^{2}(2,\delta )+1}{\eta }$, we obtain
 \begin{equation}\label{2.11}
 \begin{aligned}
   &\frac{\partial^{\alpha } }{\partial t^{\alpha }}\int _{B(x,\delta )}u^{2}dy+2\gamma \int _{B(x,\delta )}u^{2}dy\\
&\leq \Delta \int _{B(x,\delta )}u^{2}dy+2\mu (\mu C_{GN}^{2}(2,\delta )+1-\eta k)\int _{B(x,\delta )}u^{3}dy\int _{B(x,\delta )}udy\\
&\qquad+2\mu C_{GN}^{4}(2,\delta)\\
&\leq \Delta \int _{B(x,\delta )}u^{2}dy+2\mu C_{GN}^{4}(2,\delta).
\end{aligned}
 \end{equation}
 Denote $w(t)=\int_{B(x,\delta )}u^{2}dy$ the solution of the following fractional ordinary differential equation
 \begin{equation}\label{2.12}
   \begin{cases}
    _{0}^{C}\textrm{D}_{t}^{\alpha }w(t)=-2\gamma w(t)+2\mu C_{GN}^{4}(2,\delta ),\\
   w(0)=(2\delta)^{2}\left \| u_{0} \right \|_{L^{\infty}(\mathbb{R}^{N})}^{2},
   \end{cases}
 \end{equation}
  from Lemma \ref{lemma2.4}, let $A=-2\gamma,C_{3}=2\mu C_{GN}^{4}(2,\delta ),\eta =(2\delta)^{2}\left \| u_{0} \right \|_{L^{\infty}(\mathbb{R}^{N})}^{2} $,\\ for any $(x,t)\in  \mathbb{R}\times [0,T_{max})$, we obtain
 $$
 \begin{aligned}
 &\int _{B(x,\delta )}u^{2}dy=\left \| u \right \|_{L^{2}(B(x,\delta ))}^{2}\leq \left \| w(t) \right \|\\
&\leq w(0)e^{-\sigma T}+\frac{(K+\left | C_{3} \right |)T^{\alpha }e^{-\sigma T}}{\alpha \Gamma (\alpha )})Ce^{[(-2\gamma )^{\frac{1}{\alpha }}+\sigma]t}\\
&\leq (2\delta)^{2}\left \| u_{0} \right \|_{L^{\infty}(R^{N})}^{2}+\frac{(K+\left | 2\mu C_{GN}^{4}(2,\delta ) \right |)T^{\alpha }}{\alpha \Gamma (\alpha )})C:=M_{2}.
 \end{aligned}
 $$
 In conclusion, for any $(x,t)\in  \mathbb{R}\times [0,T_{max})$, we have
 \begin{equation}\label{2.13}
   \left \| u \right \|_{L^{2}(B(x,\delta ))}\leq M=\begin{cases}
	\sqrt{M_{1}}\qquad N=1\\
	\sqrt{M_{2}}\qquad N=2
\end{cases}
 \end{equation}
 and then
 \begin{equation}\label{2.14}
   \left \| u \right \|_{L^{1}B(x,\delta )}\leq (2\delta )^{\frac{N}{2}}M.
 \end{equation}
 Now we proceed to improve the $L^{2}$ boundedness of $u$ to $L^{\infty}$, which is based on the fact that for all $(x,t)\in  \mathbb{R}\times [0,T_{max})$ and Lemma \ref{lem3}-\ref{lemma2.6}, Proposition \ref{proposition.1}, we have
\begin{equation}
\begin{aligned}
  \nonumber&0\leq u(x,t)\leq \left \| \mathcal{S} _{\alpha }(t)u_{0} \right \|_{L^{\infty }(\mathbb{R}^{N})}+\int_{0}^{t}(t-s)^{\alpha -1}\left \| \mathcal{K}_{\alpha }(t-s)h(s) \right \|_{L^{\infty }(\mathbb{R}^{N})}ds\\
\nonumber&\leq C_{4}\left \| u_{0} \right \|_{L^{\infty }(\mathbb{R}^{N})}+C_{4}\int_{0}^{t}(t-s)^{\alpha -1}\left \| h(s) \right \|_{L^{\infty }(\mathbb{R}^{N})}ds\\
\nonumber&\leq C_{4}\left \| u_{0} \right \|_{L^{\infty }(\mathbb{R}^{N})}+\mu C_{4}\int_{0}^{t}(t-s)^{\alpha -1}\left \| u^{2}(s) \right \|_{L^{\infty }(\mathbb{R}^{N})}ds\\
\nonumber&\leq  C_{4}\left \| u_{0} \right \|_{L^{\infty }(\mathbb{R}^{N})}+\mu C_{4}\int_{0}^{t}(t-s)^{\alpha -1}\left \| u(s) \right \|^{2}_{L^{\infty }(\mathbb{R}^{N})}ds\\
\nonumber&\leq C_{4}\left \| u_{0} \right \|_{L^{\infty }(\mathbb{R}^{N})}+\mu C_{4}M^{2}\int_{0}^{t}(t-s)^{\alpha -1}ds\\
\nonumber&\leq C_{4}\left \| u_{0} \right \|_{L^{\infty }(\mathbb{R}^{N})}+\mu C_{4}M^{2}T^{\alpha }\frac{1}{\alpha}.
\end{aligned}
\end{equation}
So there is
\begin{equation}\label{2.22}
0\leq u(x,t)\leq C_{4}\left \| u_{0} \right \|_{L^{\infty }(\mathbb{R}^{N})}+\mu C_{4}M^{2}T^{\alpha }\frac{1}{\alpha},
\end{equation}
with
\begin{equation}
    M=\begin{cases}
	\sqrt{ 2\delta \left \| u_{0} \right \|_{L^{\infty }(\mathbb{R}^{N})}^{2}+\frac{(K+\left | 2\mu \left ( (\mu ^{\frac{1}{3}}C_{GN}^{\frac{4}{3}}(1,\delta)+1)^{6}(\eta k)^{-5}+C_{GN}^{10}(1,\delta ) \right ) \right |)T^{\alpha }}{\alpha \Gamma (\alpha )}C}, & N=1\\
	\sqrt{(2\delta)^{2}\left \| u_{0} \right \|_{L^{\infty}(\mathbb{R}^{N})}^{2}+\frac{(K+\left | 2\mu C_{GN}^{4}(2,\delta ) \right |)T^{\alpha }}{\alpha \Gamma (\alpha )}C},& N=2
\end{cases}
 \end{equation}
 defined in \eqref{2.13}. From Proposition \ref{Proposition2.1}, we obtain $T_{max}=+\infty$. The global boundedness of $u$ is obtained in time and the unique classical solution of \eqref{1}-\eqref{2} on $(x,t)\in \mathbb{R}^{N}\times [0,+\infty )$. Theorem \ref{theorem4} is thus proved.
 \begin{remark}
  This theorem is based on Chen Cheng\textsuperscript{\cite{0Global}}, citing its nonlocal terms, and processing the time fractional derivative to obtain a class of NTFRDE.In the paper, we can find that the only solution that can get equation \eqref{1}-\eqref{2} is globally bounded. And in the proof process, we use Laplace transform to solve  fractional differential equations and use the fractional Duhamel's formula (see Lemma \ref{lemma2.6}) greatly simplifies the proof of boundedness from $L^{2}$ to $L^{\infty}$.
 \end{remark}

\section{Long time behavior of solutions (the Allee effect)}

To study the long time behavior of solutions for \eqref{1}-\eqref{2}, by \cite{0Global}, we denote $$F(u):=\mu u^{2}(1-kJ*u)-\gamma u.$$ For $0<\gamma<\frac{\mu }{4k}$, there are three constant solutions for $F(u)=0$: $0,a,A$, where
\begin{equation}\label{1.7}
  a=\frac{1-\sqrt{1-4k\frac{\gamma }{\mu }}}{2k},\qquad A=\frac{1+\sqrt{1-4k\frac{\gamma }{\mu }}}{2k},
\end{equation}
 and satisfy $0<\frac{\gamma }{\mu }<a<A$.
 \begin{proposition}\label{proposition3.1}
   Under the assumptions of Theorem \ref{theorem4}, there is \\ $\left \| u(x,t) \right \|_{L^{\infty }(\mathbb{R}^{N}\times (0,+\infty ))}< a$, the function$$H(x,t)=\int _{B(x,\delta )}h(u(y,t))dy$$ with $$h(u)=Aln\left ( 1-\frac{u}{A} \right )-aln\left ( 1-\frac{u}{a} \right )$$ is nonnegative and satisfies
   \begin{equation}\label{3.1}
     \frac{\partial^{\alpha } H(x,t)}{\partial t^{\alpha }}\leq \Delta H(x,t)-D(x,t)
   \end{equation}
   with $$D(x,t)=\frac{1}{2}(A-a)\mu k\int _{B(x,\delta )}u^{2}(y,t)dy.$$
 \end{proposition}

 \begin{pro}
 Denote $k=\left \| u(x,t) \right \|_{L^{\infty }(R^{N}\times [0,\infty ))}$, then $0<k<a$. From the definition of $h(\cdot)$, it is easy to verify that
 \begin{equation}\label{3.2}
   {h}'(u)=\frac{a}{a-u}-\frac{A}{A-u}=\frac{(A-a)u}{(A-u)(a-u)},
 \end{equation}
 and
 \begin{equation}\label{3.3}
  {h}''(u)=\frac{a}{(a-u)^{2}}-\frac{A}{(A-a)^{2}}=\frac{(Aa-u^{2})(A-a)}{(A-u)^{2}(a-u)^{2}}.
 \end{equation}
 Test \eqref{1} by ${h}'(u)\varphi _{\varepsilon }$ with $\varphi _{\varepsilon }(\cdot )\in C_{0}^{\infty }(B(x,\delta )),\varphi _{\varepsilon }(\cdot )\rightarrow 1$ in $B(x,\delta)$ as $\varepsilon \rightarrow 0$. Integrating by parts over $B(x,\delta)$ and $\Delta h(u)={h}''(u)\left | \bigtriangledown u \right |^{2}+{h}'(u)\Delta u$, we obtain
 $$
 \begin{aligned}
& \int _{B(x,\delta )}{h}'(u)\varphi _{\varepsilon }\partial _{t}^{\alpha }udy\\
&=\int _{B(x,\delta )}\Delta u{h}'(u)\varphi _{\varepsilon }dy+\int _{B(x,\delta )}[\mu u^{2}(1-kJ*u)-\gamma u]{h}'(u)\varphi _{\varepsilon }dy\\
&= \int _{B(x,\delta )}\Delta h(u)\varphi _{\varepsilon }dy-\int _{B(x,\delta )}{h}''(u)\left | \bigtriangledown u \right |^{2}\varphi _{\varepsilon }dy\\
&\qquad+\int _{B(x,\delta )}[\mu u^{2}(1-kJ*u)-\gamma u]{h}'(u)\varphi _{\varepsilon }dy.
 \end{aligned}
 $$
 From Theorem \ref{theorem2} and Equation \eqref{4}, we get
 $$\int _{B(x,\delta )}{h}'(u)\partial _{t}^{\alpha }udy\geq \int _{B(x,\delta )}\partial _{t}^{\alpha }h(u)dy=\partial _{t}^{\alpha }\int _{B(x,\delta )}h(u)dy.$$
 And by Lemma \ref{lemma2.8}, noticing
 $$
 \begin{aligned}
 \int _{B(x,\delta )}\Delta h(u)dy&=\int _{B(0,\delta )}\Delta h(u(y+x),t)dy \\
 &=\Delta \int _{B(0,\delta )}h(u(x+y),t)dy=\Delta \int _{B(x,\delta )}h(u)dy,
 \end{aligned}
 $$
 taking $\varepsilon \rightarrow 0$, we obtain
 $$
 \begin{aligned}
& \frac{\partial^{\alpha } }{\partial t^{\alpha }} \int _{B(x,\delta )}h(u)dy\\
&\leq \Delta \int _{B(x,\delta )}h(u)dy-\int _{B(x,\delta )}{h}''(u)\left | \bigtriangledown u \right |^{2}dy\\
&\qquad+\int _{B(x,\delta )}[\mu u^{2}(1-kJ*u)-\gamma u]{h}'(u)dy,
 \end{aligned}
 $$
 which is
 \begin{equation}\label{3.4}
 \begin{aligned}
   &\frac{\partial^{\alpha } }{\partial t^{\alpha }}H(x,t)\\
 &\leq \Delta H(x,t)-\int _{B(x,\delta )}{h}''(u)\left | \bigtriangledown u \right |^{2}dy\\
 &\qquad+\int _{B(x,\delta )}[\mu u^{2}(1-kJ*u)-\gamma u]{h}'(u)dy.
\end{aligned}
 \end{equation}
 By \eqref{1.7}, we can get $$\mu u^{2}(1-ku)-\gamma u=k\mu u(A-u)(u-a),$$
and
\begin{equation}\label{3.5}
\begin{aligned}
&\int _{B(x,\delta )}{h}'(u)[\mu u^{2}(1-kJ*u)-\gamma u]dy\\
&=\int _{B(x,\delta )}{h}'(u)[\mu u^{2}(1-ku)-\gamma u]dy+\mu k\int _{B(x,\delta )}{h}'(u)u^{2}(u-J*u)dy\\
&=-(A-a)\mu k\int _{B(x,\delta )}u^{2}(y,t)dy+\mu k\int _{B(x,\delta )}{h}'(u)u^{2}(u-J*u)dy.
\end{aligned}
\end{equation}
Noticing that when $0\leq u\leq k$, there is $$0\leq {h}'(u)u\leq \frac{(A-a)k^{2}}{(A-k)(a-k)}.$$
From Young's inequality and the median value theorem, we can get
\begin{equation}\label{3.6}
\begin{aligned}
&\mu k\int _{B(x,\delta )}{h}'(u)u^{2}(u-J*u)dy\\
&\leq \mu k\int _{B(x,\delta )}\int _{B(x,\delta )}{h}'(u)u^{2}(y,t)(u(y,t)-u(z,t))J(y-z)dzdy\\
&\leq \frac{(A-a)K^{2}}{(A-K)(a-K)}\mu k\int _{B(x,\delta )}\int _{B(x,\delta )}u(y,t)\left | (u(y,t)-u(z,t)) \right |J(y,z)dzdy\\
&\leq \frac{(A-a)K^{4}\mu k}{2(A-K)^{2}(a-K)^{2}}\int _{B(x,\delta )}\int _{B(x,\delta )}(u(z,t)-u(y,t))^{2}J(z-y)dzdy\\
&\qquad+\frac{1}{2}(A-a)\mu k\int _{B(x,\delta )}u^{2}(y,t)dy\\
&\leq \frac{(A-a)K^{4}\mu k}{2(A-K)^{2}(a-K)^{2}}\int _{B(x,\delta )}\int _{B(x,\delta )}\int_{0}^{1}\left | \bigtriangledown u(y+\theta (z-y),t) \right |^{2}\left | z-y \right |^{2}\\
&\qquad J(z-y) d\theta dzdy+\frac{1}{2}(A-a)\mu k\int _{B(x,\delta )}u^{2}(y,t)dy,
  \end{aligned}
\end{equation}
changing the variables ${y}'=y+\theta (z-y) , {z}'=z-y$, then$$\begin{vmatrix}
	& \frac{\partial{y}' }{\partial y} \qquad \frac{\partial {y}'}{\partial z}\\
	& \frac{\partial {z}'}{\partial y}\qquad  \frac{\partial {z}'}{\partial z}\\
\end{vmatrix}=\begin{vmatrix}
	1-\theta   \qquad\theta & \\
	-1\qquad	1 &
\end{vmatrix}=1-\theta+\theta=1.$$
For any $\theta \in  [0,1],y,z\in  B(x,\delta)$, we have ${y}'\in B((1-\theta )x+\theta z,(1-\theta )\delta ) ,{z}'\in B(x-y,\delta )$. Noticing $B((1-\theta )x+\theta z,(1-\theta )\delta )\subseteq B(x,\delta)$ and $B(x-y,\delta)\subseteq B(0,2\delta)$, we obtain
\begin{equation}\label{3.7}
\begin{aligned}
&\frac{(A-a)K^{4}\mu k}{2(A-K)^{2}(a-K)^{2}}\int _{B(x,\delta )}\int _{B(x,\delta )}\int_{0}^{1}\left | \bigtriangledown u(y+\theta (z-y),t) \right |^{2}\left | z-y \right |^{2}\\
& \qquad J(z-y)d\theta dzdy\\
&\leq \frac{(A-a)K^{4}\mu k}{2(A-K)^{2}(a-K)^{2}}\int_{0}^{1}d\theta \int _{B(x,\delta )}\int _{B(x,\delta )}\left | \bigtriangledown u({y}',t) \right |^{2}\left | {z}' \right |^{2}J({z}')d{z}'d{y}'\\
&\leq \frac{(A-a)K^{4}\mu k (2\delta )^{2}}{2(A-K)^{2}(a-K)^{2}}\int _{B(x,\delta )}\left | \bigtriangledown u(y,t) \right |^{2}dy.
\end{aligned}
\end{equation}
Combining \eqref{3.5}-\eqref{3.7}, we obtain
\begin{equation}\label{3.8}
\begin{aligned}
  &\int _{B(x,\delta )}{h}'(u)[\mu u^{2}(1-kJ*u)-\gamma u]dy\\
&\leq -\frac{1}{2}(A-a)\mu k\int _{B(x,\delta )}u^{2}(y,t)dy\\
&\qquad+\frac{(A-a)k^{4}\mu k(2\delta )^{2}}{2(A-a)^{2}(a-k)^{2}}\int _{B(x,\delta )}\left | \bigtriangledown u(y,t) \right |^{2}dy.
\end{aligned}
\end{equation}
From \eqref{3.3}, noticing $0\leq u<a$, we obtain$${h}''(u)> \frac{(A-a)^{2}}{A^{2}a},$$
inserting \eqref{3.8} into \eqref{3.4}, we obtain
\begin{equation}\label{3.9}
  \begin{aligned}
  \frac{\partial^{\alpha } }{\partial t^{\alpha }}H(x,t)&\leq \Delta H(x,t)+[-\frac{(A-a)^{2}}{A^{2}a}+\frac{(A-a)k^{4}\mu k(2\delta )^{2}}{2(A-k)^{2}(a-k)^{2}}]\\
  &\qquad \int _{B(x,\delta )}\left | \bigtriangledown u(y,t) \right |^{2}dy-\frac{1}{2}(A-a)\mu k\int _{B(x,\delta )}u^{2}(y,t)dy.
  \end{aligned}
\end{equation}
By choosing $\delta$ sufficiently small such that$$-\frac{(A-a)^{2}}{A^{2}a}+\frac{(A-a)k^{4}\mu k(2\delta )^{2}}{2(A-k)^{2}(a-k)^{2}}\leq 0,$$ then
$$\frac{\partial^{\alpha } }{\partial t^{\alpha }}H(x,t)\leq \Delta H(x,t)-\frac{1}{2}(A-a)\mu k\int _{B(x,\delta )}u^{2}(y,t)dy,$$
making $$D(x,t)=\frac{1}{2}(A-a)\mu k\int _{B(x,\delta )}u^{2}(y,t)dy.$$
 \end{pro}
\begin{remark}
  The proof of the inequality is based on the proof method in \cite{0Global}, in which the scaling and construction of nonlocal terms are complete. Therefore, the paper uses Theorem \ref{theorem2} and Lemma \ref{lemma23}, and Median theorem to obtain fractional differential inequality on the basis of changing its conditions. This Proposition is particularly important in the later proof process.
\end{remark}

\begin{definition}\textsuperscript{\cite{2006Global}}
  Let $X$ be a Banach space,$z_{0}$ belong to $X$, and $f\in L^{1}(0,T;X)$. The function $z(x,t)\in C([0,T];X)$ given by
  \begin{equation}\label{pp}
      z(t)=e^{-t}e^{t\Delta}z_{0}+\int_{0}^{t}e^{-(t-s)}\cdot e^{(t-s)\Delta }f(s)ds,\qquad 0\leq t\leq T,
  \end{equation}
  is the mild solution of \eqref{pp} on $[0,T]$, where $(e^{t\Delta}f)(x,t)=\int _{\mathbb{R}^{N}}G(x-y,t)f(y)dy$ and $G(x,t)$ is the heat kernel by $G(x,t)=\frac{1}{(4\pi t)^{N/2}}exp(-\frac{\left | x \right |^{2}}{4t})$.
\end{definition}
\begin{lemma}\textsuperscript{\cite{2006Global}}\label{lemma5}
  Let $0\leq q\leq p\leq \infty ,\frac{1}{q}-\frac{1}{p}<\frac{1}{N}$ and suppose that $z$ is the function given by \eqref{pp} and $z_{0}\in W^{1,p}(\mathbb{R}^{N})$. If $f\in L^{\infty }(0,\infty ;L^{q}(\mathbb{R}^{N}))$, then
  \begin{eqnarray}
   \nonumber \left \| z(t) \right \|_{L^{p}(\mathbb{R}^{N})}&\leq& \left \|  z_{0}\right \|_{L^{p}(\mathbb{R}^{N})}+C\cdot \Gamma (\gamma )\underset{0<s<t}{sup}\left \| f(s) \right \|_{L^{q}(\mathbb{R}^{N})},\\
\nonumber\left \| \bigtriangledown z(t) \right \|_{L^{p}(\mathbb{R}^{N})}&\leq& \left \|  \bigtriangledown z_{0}\right \|_{L^{p}(\mathbb{R}^{N})}+C\cdot \Gamma (\tilde{\gamma })\underset{0<s<t}{sup}\left \| f(s) \right \|_{L^{q}(\mathbb{R}^{N})},
  \end{eqnarray}
  for $t\in [0,\infty )$, where $C$ is a positive constant independent of $p,\Gamma (\cdot )$ is the gamma function, and $\gamma =1-(\frac{1}{q}-\frac{1}{p})\cdot \frac{N}{2},\tilde{\gamma }=\frac{1}{2}-(\frac{1}{q}-\frac{1}{p})\cdot \frac{N}{2}$.
\end{lemma}

 \noindent\textbf{Proof of Theorem \ref{theorem5}}. From the proof of Theorem \ref{theorem4}, for any $\forall{K}'>0$ , from \eqref{2.22} and the definition of
 $M$, there exist $\mu ^{*}({K}')>0$ and $m_{0}({K}')>0$ such that for $\mu \in  (0,\mu ^{*}({K}'))$ and $\left \| u_{0} \right \|_{L^{\infty }(\mathbb{R}^{N})}<m_{0}({K}')$, then for any $\left ( x,t \right )\in \mathbb{R}^{N}\times [0,+\infty ),M$ is sufficiently small such that $$\left \| u(x,t) \right \|_{L^{\infty }(\mathbb{R}^{N}\times [0,+\infty ))}< {K}'.$$
 (i) \textbf{The case}\, $\left \| u(x,t) \right \|_{L^{\infty }(\mathbb{R}^{N}\times [0,+\infty ))}<{K}'=\frac{\gamma }{\mu }$.\\
 Noticing $\sigma =\gamma -\mu \left \| u(x,t) \right \|_{L^{\infty }(\mathbb{R}^{N}\times [0,\infty ))}>0$, we have
 \begin{equation}\label{3.10}
   \begin{aligned}
   &\frac{\partial^{\alpha } u}{\partial t^{\alpha }}=\Delta u+\mu u^{2}(1-kJ*u)-\gamma u\\
&=\Delta u+u[\mu u(1-kJ*u)-\gamma ]\\
&=\Delta u-u[\gamma -\mu u(1-kJ*u)]\\
&\leq \Delta u-u[\gamma -\mu u]\\
&\leq \Delta u-u\sigma.
   \end{aligned}
 \end{equation}
 Consider the following fractional ordinary differential equation
 \begin{equation}
   \begin{cases}
   _{0}^{C}\textrm{D}_{t}^{\alpha }w(t)=-\sigma w,\\
	w(0)=\left \| u_{0} \right \|_{L^{\infty }(\mathbb{R}^{N})}.
   \end{cases}
 \end{equation}
 By Lemma \ref{lemma2.4} and Equation \eqref{2.6}, we obtain
 $$u(x,t)\leq w(t)=w(0)E_{\alpha ,1}(-\sigma t)\leq \left \| u_{0} \right \|_{L^{\infty }(\mathbb{R}^{N})}e^{(-\sigma)^{\frac{1}{\alpha }} t}$$
 for all $(x,t)\in \mathbb{R}^{N}\times [0,+\infty )$ and obtain $$\left \| u(\cdot ,t) \right \|_{L^{\infty }(\mathbb{R}^{N})}\leq \left \| u_{0} \right \|_{L^{\infty }(\mathbb{R}^{N})}e^{(-\sigma)^{\frac{1}{\alpha }} t}.$$

 (ii)\textbf{The case}\, $\left \| u(x,t) \right \|_{L^{\infty }(\mathbb{R}^{N}\times [0,+\infty ))}< {K}'=a$.\\
 Denoting $H_{0}(y)=H(y,0)$, from \eqref{3.1} in Proposition \ref{proposition3.1} and fractional Duhamel’s formula \ref{lemma2.6}, for any $(x,t)\in \mathbb{R}^{N}\times [0,+\infty )$, we have
$$ H(x,t)\leq \left \| H_{0} \right \|_{L^{\infty }(\mathbb{R}^{N})}-\int_{0}^{t}(t-s)^{\alpha -1}\mathcal{K}_{\alpha }(t-s)D(y,s)ds,$$
from which we obtain$$\int_{0}^{t}(t-s)^{\alpha -1}\mathcal{K}_{\alpha }(t-s)D(y,s)ds \leq \left \| H_{0} \right \|_{L^{\infty }(\mathbb{R}^{N})}.$$
Due to the fact that $u$ is a classical solution, we have that $$\mathcal{K}_{\alpha }(t-s)D(y,s)\in C^{2,1}(\mathbb{R}^{N}\times [0,\infty )),$$
which implies that for all $x \in \mathbb{R}^{N}$, the following limit holds:
$$\lim_{t\rightarrow \infty }\lim_{s\rightarrow t}\mathcal{K}_{\alpha }(t-s)D(y,s)=0,$$
or equivalently
$$\lim_{t\rightarrow \infty }\lim_{s\rightarrow t}\mathcal{K}_{\alpha }(t-s)\int _{B(x,\delta )}u^{2}(z,s)dz=0,$$
which together with the fact that the heat kernel converges to delta function as $s\rightarrow t$, we have that for any $x \in \mathbb{R}^{N}$,
$$\lim_{t\rightarrow \infty }\int _{B(x,\delta )}u^{2}(y,t)dy=0.$$
Furthermore, with the uniform boundedness of $u$ on $\mathbb{R}^{N}\times[0,\infty)$, following Lemma \ref{lemma5}, we can obtain the global boundedness of
$\left \| \bigtriangledown u(\cdot ,t) \right \|_{L^{\infty }(\mathbb{R}^{N})}$, from which and Lemma \ref{lemma2.1} with $\Omega =B(x,\delta ),p=q=\infty ,r=2$, the convergence of $\left \|  u(\cdot ,t)\right \|_{L^{\infty }(B(x,\delta ))}$ follows from the convergences of $\left \|  u(\cdot ,t)\right \|_{L^{2}(B(x,\delta ))}$ immediately. This is, we obtain $$\left \|  u(\cdot ,t)\right \|_{L^{\infty }(B(x,\delta ))}\rightarrow 0$$ as $t \rightarrow 0$.
Therefore, for any compact set in $\mathbb{R}^{N}$, by finite covering, we obtain that $u$ converges to $0$ uniformly in that compact set, which means that $u$ converges locally uniformly to $0$ in $\mathbb{R}^{N}$ as $\rightarrow \infty$. The proof is complete.

 \section{Global bounded of solutions for a nonlinear NTFRDE}

 \begin{lemma}\textsuperscript{\cite{2016Weak}}\label{lemma2.14}
 Suppose that a nonnegative function $u(t)\geq 0$  satisfies
 \begin{equation}\label{44}
   _{0}^{C}\textrm{D}_{t}^{\alpha }u(t)+c_{1}u(t)\leq f(t)
 \end{equation}
  for almost all $t\in [0,T]$, where $c_{1}>0$, and the function $f(t)$ is nonnegative and integrable for $t\in [0,T]$. Then
  \begin{equation}\label{55}
    u(t)\leq u(0)+\frac{1}{\Gamma (\alpha )}\int_{0}^{t}(t-s)^{\alpha -1}f(s)ds.
  \end{equation}
 \end{lemma}

\begin{corollary}\label{corollary1}
  Suppose that a nonnegative function $u(t)\geq 0$  satisfies
  \begin{equation}
    \nonumber _{0}^{C}\textrm{D}_{t}^{\alpha }u(t)+c_{1}u(t)\leq b
  \end{equation}
for almost all $t\in [0,T]$, where $b,c_{1}>0$ are all contants. then
\begin{equation}\label{585}
   u(t)\leq u(0)+\frac{bT^{\alpha }}{\alpha \Gamma (\alpha )}.
  \end{equation}
\end{corollary}

\begin{pro}
  Let $ _{0}^{C}\textrm{D}_{t}^{\alpha }u(t)+c_{1}u(t)=g(t)$, then by \cite{Kilbas2006}, we have the following expression:
  $$v(t)=v(0)E_{\alpha }(-c_{1}t^{\alpha })+\int_{0}^{t}(t-s)^{\alpha -1}E_{\alpha ,\alpha }(-C_{1}(t-s)^{\alpha })g(s)ds,$$
  where$$E_{\alpha }(t)=\sum_{k=0}^{\infty }\frac{t^{k}}{\Gamma (\alpha k+1)},\qquad E_{\alpha,\alpha }(t)=\sum_{k=0}^{\infty }\frac{t^{k}}{\Gamma (\alpha k+\alpha)},$$ are Mittag-Leffler functions. By Lemma \ref{lemma288}-\ref{lemma299}, we know
  $0<E_{\alpha }(-c_{1}t^{\alpha })\leq 1$,\quad$0<E_{\alpha ,\alpha }(-C_{1}(t-s)^{\alpha })\leq \frac{1}{\Gamma (\alpha )}$ and $g(t)\leq f(t)$, then
  \begin{equation}
    \begin{aligned}
  \nonumber  v(t)&\leq v(0)E_{\alpha }(-c_{1}t^{\alpha })+\frac{1}{\Gamma (\alpha )}\int_{0}^{t}(t-s)^{\alpha -1}g(s)ds\\
\nonumber  &\leq v(0)+\frac{1}{\Gamma (\alpha )}\int_{0}^{t}b(t-s)^{\alpha -1}ds=u(0)+\frac{bT^{\alpha }}{\alpha \Gamma (\alpha )}.
    \end{aligned}
  \end{equation}
\end{pro}
The proof of this corollary has been proved.

\begin{remark}
  Fractional differential inequalities are commonly used inequalities in fractional differential equations, and the proof method is mentioned in \cite{2016Weak}. But compared with Lemma \ref{lemma2.14}, it replaces the function term with a special form of constant term, which leads to inferences.
\end{remark}

 \begin{lemma}\label{lemma444}
   Assume the function $y_{k}(t)$ is nonnegative and exists the Caputo fractional derivative for $t\in [0,T]$ satisfying
   \begin{equation}\label{4.17}
  _{0}^{C}\textrm{D}_{t}^{\alpha }y_{k}(t)\leq -y_{k}+a_{k}(y_{k-1}^{\gamma _{1}}(t)+y_{k-1}^{\gamma _{2}}(t)),
   \end{equation}
   where $a_{k}=\bar{a}3^{rk}>1$ with $\bar{a},r$ are positive bounded constants and $0<\gamma _{2}<\gamma _{1}\leq 3$. Assume also that there exists a bounded constant $K\geq 1$ such that $y_{k}(0)\leq K^{3^{k}}$, then
   \begin{equation}\label{4.18}
    y_{k}(t)\leq (2\bar{a})^{\frac{3^{k}-1}{2}}3^{r (\frac{3^{k+1}}{4}-\frac{k}{2}-\frac{3}{4})}max\left \{ \underset{t\in [0,T]}{sup} y_{0}^{3^{k}}(t),K^{3^{k}}\right \}\frac{T^{\alpha }}{\alpha \Gamma (\alpha)}.
   \end{equation}
 \end{lemma}

 \begin{pro}
 From Lemma \ref{lemma2.14} and Eqution \eqref{44}-\eqref{585}, we make $$ c_{1}=-1,f(t)=a_{k}(y_{k-1}^{\gamma _{1}}(t)+y_{k-1}^{\gamma _{2}}(t))\leq 2a_{k}max\left \{ 1,\underset{t\in [0,T]}{sup} y_{k-1}^{3}(t)\right \},$$
  then
 \begin{equation}\label{4.19}
 \begin{aligned}
& y_{k}(t)\leq y_{k}(0)+\frac{1}{\Gamma (\alpha )}\int_{0}^{t}(t-s)^{\alpha -1}f(s)ds\\
&\leq K^{3^{k}}+\frac{2a_{k}}{\Gamma (\alpha )}\int_{0}^{t}(t-s)^{\alpha -1}max\left \{ 1,\underset{t\in [0,T]}{sup} y_{k-1}^{3}(s)\right \}ds\\
&\leq K^{3^{k}}+\frac{2a_{k}}{\Gamma (\alpha )}max\left \{ 1,\underset{t\in [0,T]}{sup} y_{k-1}^{3}(t)\right\}\int_{0}^{t}(t-s)^{\alpha -1}ds\\
&\leq K^{3^{k}}+\frac{T^{\alpha }}{\alpha \Gamma (\alpha )}2a_{k}max\left \{ 1,\underset{t\in [0,T]}{sup} y_{k-1}^{3}(t)\right\}\\
&\leq \frac{T^{\alpha }}{\alpha \Gamma (\alpha )}2a_{k}max\left \{K^{3^{k}},\underset{t\in [0,T]}{sup} y_{k-1}^{3}(t)\right\}.
 \end{aligned}
 \end{equation}
 Then from \eqref{4.19} after some iterative steps we have
 \begin{equation}
   \begin{aligned}
   \nonumber &y_{k}(t)\leq 2a_{k}(2a_{k-1})^{3}(2a_{k-2})^{3^{2}}(2a_{k-3})^{3^{3}}\cdots (2a_{1})^{3^{k-1}}max\left \{K^{3^{k}},\underset{t\in [0,T]}{sup} y_{0}^{3^{k}}(t)\right\}\\
  & \qquad \frac{T^{\alpha }}{\alpha \Gamma (\alpha )}\\
   \nonumber&=(2\bar{a})^{1+3+3^{2}+3^{3}+\cdots 3^{k-1}}3^{r(k+3(k-1)+3^{2}(k-2)+\cdots +3^{k-1})}max\left \{K^{3^{k}},\underset{t\in [0,T]}{sup} y_{0}^{3^{k}}(t)\right\}\\
   &\qquad \frac{T^{\alpha }}{\alpha \Gamma (\alpha )}\\
   \nonumber&=(2\bar{a})^{\frac{3^{k}-1}{2}}3^{r(\frac{3^{k+1}}{4}-\frac{k}{2}-\frac{3}{4})}max\left \{K^{3^{k}},\underset{t\in [0,T]}{sup} y_{0}^{3^{k}}(t)\right\}\frac{T^{\alpha }}{\alpha \Gamma (\alpha )}.
   \end{aligned}
 \end{equation}
 \end{pro}

 \begin{lemma}\textsuperscript{\cite{1959On}}\textsuperscript{\cite{1971Propriet}}\label{lemma33}
   When the parameters $p,q,r$ meet any of the following conditions:

   (i) $q >N \geq 1$, $r\geq1$ and $p=\infty$;

   (ii)$q>max\left \{ 1,\frac{2N}{N+2} \right \},1\leq r<\sigma $ and $r<p<\sigma +1$ in
   \begin{equation}
    \nonumber \sigma :=\begin{cases}
               \frac{(q-1)N+q}{N-q},&q<N,\\
               \infty ,&q\geq N.
              \end{cases}
   \end{equation}
Then the following inequality is established
$$\left \| u \right \|_{L^{p}(\mathbb{R}^{N})}\leq C_{GN}\left \| u \right \|_{L^{r}(\mathbb{R}^{N})}^{1-\lambda ^{*}}
\left \| \bigtriangledown u \right \|_{L^{q}(\mathbb{R}^{N})}^{\lambda ^{*}}$$
among
$$\lambda ^{*}=\frac{qN(p-r)}{p[N(q-r)+qr]}.$$
 \end{lemma}

\begin{lemma}\textsuperscript{\cite{Shen2016A}}\label{lemma333}
  Let $N\geq 1$. $p$ is the exponent from the Sobolev embedding theorem, $i.e.$
  \begin{equation}\label{22.1}
    \begin{cases}
      p=\frac{2N}{N-2},& N\geq 3\\
     2<p<\infty,& N=2\\
    p=\infty,& N=1
    \end{cases}
  \end{equation}
  $1\leq r<q<p$ and $\frac{q}{r}< \frac{2}{r}+1-\frac{2}{p}$, then for $v\in {H}'(\mathbb{R}^{N})$ and $v\in L^{r}(\mathbb{R}^{N})$, it holds
  \begin{equation}\label{99.1}
  \begin{aligned}
    \left \| v \right \|_{L^{q}(\mathbb{R}^{N})}^{q}&\leq
     C(N)c^{\frac{\lambda q}{2-\lambda q}}_{0}\left \| v \right \|_{L^{p}(\mathbb{R}^{N})}^{\gamma }+c_{0}\left \| \bigtriangledown u \right \|_{L^{2}(\mathbb{R}^{N})}^{2},\quad N>2,\\
\left \| v \right \|_{L^{q}(\mathbb{R}^{N})}^{q}&\leq C(N)(c^{\frac{\lambda q}{2-\lambda q}}_{0}+c^{-\frac{\lambda q}{2-\lambda q}}_{1})\left \| v \right \|_{L^{r}(\mathbb{R}^{N})}^{\gamma }+c_{0}\left \| \bigtriangledown u \right \|_{L^{2}(\mathbb{R}^{N})}^{2}\\
&\qquad+c_{1}\left \| v \right \|_{L^{2}(\mathbb{R}^{N})}^{2},\quad  N=1,2.
	\end{aligned}
  \end{equation}
  Here $C(N)$ are constants depending on $N,c_{0},c_{1}$ are arbitrary positive constants and
  \begin{equation}\label{112}
  \nonumber \lambda =\frac{\frac{1}{r}-\frac{1}{q}}{\frac{1}{r}-\frac{1}{p}}\in (0,1),\gamma =\frac{2(1-\lambda )q}{2-\lambda q}=\frac{2(1-\frac{q}{p})}{\frac{2-q}{r}-\frac{2}{p}+1}.
  \end{equation}
\end{lemma}

\noindent\textbf{Proof of Theorem \ref{theorem6}}

\textbf{Step 1. $L^{p}$ estimates}:

\noindent For any $x\in \mathbb{R}^{N}$, multiply \eqref{1.1.4} by $u^{k-1},k>1$ and integrating by parts over $\mathbb{R}^{N}$. From Theorem \ref{theorem1} and Corollary \ref{corollary2}, we have$$u^{k-1}(t)\int _{\mathbb{R}^{N}}(_{0}^{C}\textrm{D} _{t}^{\alpha }u)dx\geq \frac{1}{k}(_{0}^{C}\textrm{D} _{t}^{\alpha }\int _{\mathbb{R}^{N}}u^{k}dx).$$
Thus, we obtain
$$
\begin{aligned}
\frac{1}{k}(_{0}^{C}\textrm{D} _{t}^{\alpha }\int _{\mathbb{R}^{N}}u^{k}dx)&\leq-\int _{\mathbb{R}^{N}}(u+1)^{m-1}\bigtriangledown u\bigtriangledown (u^{k-1})dx\\
&\qquad+\int _{\mathbb{R}^{N}}u^{k+1}dx(1-\int _{\mathbb{R}^{N}}udx)
-\int _{\mathbb{R}^{N}}u^{k}dx\\
&\leq-(k-1)\int _{\mathbb{R}^{N}}(u+1)^{m-1}u^{k-2}\left | \bigtriangledown u \right |^{2}dx\\
&\qquad+\int _{\mathbb{R}^{N}}u^{k+1}dx(1-\int _{\mathbb{R}^{N}}udx)-\int _{\mathbb{R}^{N}}u^{k}dx\\
&\leq -(k-1)\int _{\mathbb{R}^{N}}u^{m+k-3}\left | \bigtriangledown u \right |^{2}dx\\
&\qquad+\int _{\mathbb{R}^{N}}u^{k+1}dx(1-\int _{\mathbb{R}^{N}}udx)
-\int _{\mathbb{R}^{N}}u^{k}dx.
\end{aligned}
$$
Present proof $$\int _{R^{N}}u^{m+k-3}\left | \bigtriangledown u \right |^{2}dx=\frac{4}{(m+k-1)^{2}}\left \| \bigtriangledown u^{\frac{m+k-1}{2}} \right \|_{L^{2}(\mathbb{R}^{N})}^{2}.$$
From
$$
\begin{aligned}
\left \| \bigtriangledown u^{\frac{m+k-1}{2}} \right \|_{L^{2}(\mathbb{R}^{N})}^{2}&=\int _{\mathbb{R}^{N}}(\bigtriangledown u^{\frac{m+k-1}{2}})^{2}dx\\
&=\int _{\mathbb{R}^{N}}\bigtriangledown u^{\frac{m+k-1}{2}}\cdot \bigtriangledown u^{\frac{m+k-1}{2}}dx\\
&=\int _{\mathbb{R}^{N}}\frac{(m+k-1)^{2}}{4}u^{m+k-3}\left | \bigtriangledown u \right |^{2},
\end{aligned}
$$
then the following inequality holds
$$
\begin{aligned}
\frac{1}{k}(_{0}^{C}\textrm{D} _{t}^{\alpha }\int _{\mathbb{R}^{N}}u^{k}dx)&\leq-\frac{4(k-1)}{(m+k-1)^{2}}\left \| \bigtriangledown u^{\frac{m+k-1}{2}} \right \|_{L^{2}(\mathbb{R}^{N})}^{2}+\int _{\mathbb{R}^{N}}u^{k+1}dx\\
&\quad -\int _{\mathbb{R}^{N}}u^{k+1}dx\int _{R^{N}}udx-\int _{\mathbb{R}^{N}}u^{k}dx.
\end{aligned}
$$
Move the term to simplify and get
\begin{equation}\label{4.1.1}
  \begin{aligned}
  &_{0}^{C}\textrm{D} _{t}^{\alpha }\int _{\mathbb{R}^{N}}u^{k}dx+k\frac{4(k-1)}{(m+k-1)^{2}}\left \| \bigtriangledown u^{\frac{m+k-1}{2}} \right \|_{L^{2}(\mathbb{R}^{N})}^{2}\\
  &\qquad+k\int _{\mathbb{R}^{N}}u^{k+1}dx\int _{\mathbb{R}^{N}}udx+k\int _{\mathbb{R}^{N}}u^{k}dx\\
&\leq k\int _{\mathbb{R}^{N}}u^{k+1}dx.
  \end{aligned}
\end{equation}
The following estimate $k\int _{\mathbb{R}^{N}}u^{k+1}dx$. when $m\leq3$ and
\begin{equation}\label{5.5}
 k> max\left \{ \frac{1}{4}(3-m)(N-2) -(m-1),(1-\frac{m}{2})N-1,N(2-m)-2\right \},
\end{equation}
from Lemma \ref{lemma33}, we obtain
$$
\begin{aligned}
&k\int _{\mathbb{R}^{N}}u^{k+1}dx=k\left \| u^{\frac{k+m-1}{2}} \right \|_{L^{\frac{2(k+1)}{k+m+1}}(\mathbb{R}^{N})}^{\frac{2(k+1)}{k+m-1}}\\
&\leq k\left \| u^{\frac{k+m-1}{2}} \right \|_{L^{\frac{2(k+1)}{k+m+1}}(\mathbb{R}^{N})}^{\frac{2(k+1)}{k+m-1}\lambda ^{*}}\left \| \bigtriangledown u^{\frac{k+m-1}{2}} \right \|_{L^{2}(\mathbb{R}^{N})}^{\frac{2(k+1)(1-\lambda ^{*})}{k+m-1}}.
\end{aligned}
$$
Knowing from \eqref{5.5}
$$\lambda ^{*}=\frac{1+\frac{(k+m-1)N}{2(k+1)}-\frac{N}{2}}{1+\frac{(k+m-1)N}{k+2}-\frac{N}{2}}\epsilon \left \{ max\left \{ 0,\frac{2-m}{k+1} \right \},1 \right \},$$
using Young's inequality, there are
\begin{equation}\label{4.1.3}
\begin{aligned}
k\int _{\mathbb{R}^{N}}u^{k+1}dx&\leq k\left \| u^{\frac{k+m-1}{2}} \right \|_{L^{\frac{k+2}{k+m+1}}(\mathbb{R}^{N})}^{\frac{2(k+1)}{k+m-1}\lambda ^{*}}\left \| \bigtriangledown u^{\frac{k+m-1}{2}} \right \|_{L^{2}(R^{N})}^{\frac{2(k+1)(1-\lambda ^{*})}{k+m-1}}\\
&\leq \frac{2k(k-1)}{(m+k-1)^{2}}\left \| \bigtriangleup u^{\frac{k+m-1}{2}} \right \|_{L^{2}(R^{N})}^{2}\\
&\qquad+C_{1}(N,k,m)\left \| u^{\frac{k+m-1}{2}} \right \|_{L^{\frac{k+2}{k+m-1}}(\mathbb{R}^{N})}^{Q}
\end{aligned}
\end{equation}
where is $Q=\frac{2(k+1)\lambda ^{*}}{m-2+(k+1)\lambda ^{*}}$.
Next estimate $\left \| u^{\frac{k+m-1}{2}} \right \|_{L^{\frac{k+2}{k+m-1}}(\mathbb{R}^{N})}^{Q}$.
We will use the interpolation inequality to get
\begin{equation}\label{4.1.4}
  \begin{aligned}
  \left \| u^{\frac{k+m-1}{2}} \right \|_{L^{\frac{k+2}{k+m-1}}(\mathbb{R}^{N})}^{Q}&\leq\left \| u^{\frac{k+m-1}{2}} \right \|_{L^{\frac{2(k+1)}{k+m-1}}(\mathbb{R}^{N})}^{Q\lambda }\left \| u^{\frac{k+m-1}{2}} \right \|_{L^{\frac{2}{k+m-1}}(\mathbb{R}^{N})}^{Q(1-\lambda )}\\
&\leq(\left \| u^{\frac{k+m-1}{2}} \right \|_{L^{\frac{2(k+1)}{k+m-1}}(\mathbb{R}^{N})}^{\frac{2(k+1)}{k+m-1} }\left \| u^{\frac{k+m-1}{2}} \right \|_{L^{\frac{2}{k+m-1}}(\mathbb{R}^{N})}^{\frac{2}{k+m-1}})^{\frac{Q\lambda (k+m-1)}{2(k+1)}}\\
&\left \| u^{\frac{k+m-1}{2}} \right \|_{L^{\frac{2}{k+m-1}}(\mathbb{R}^{N})}^{Q(1-\lambda -\frac{\lambda }{k+1})}
  \end{aligned}
\end{equation}
in $\lambda =\frac{k+1}{k+2}$, and $$Q(1-\lambda -\frac{\lambda }{k+1})=0.$$ Then
\begin{equation}\label{4.1.5}
\begin{aligned}
  &C_{1}(N,k,m)\left \| u^{\frac{k+m-1}{2}} \right \|_{L^{\frac{k+2}{k+m-1}}(\mathbb{R}^{N})}^{Q}\\
  &\leq C_{1}(N,k,m)(\left \| u^{\frac{k+m-1}{2}} \right \|_{L^{\frac{2(k+1)}{k+m-1}}(\mathbb{R}^{N})}^{\frac{2(k+1)}{k+m-1} }\left \| u^{\frac{k+m-1}{2}} \right \|_{L^{\frac{2}{k+m-1}}(R^{N})}^{\frac{2}{k+m-1}})^{\frac{Q\lambda (k+m-1)}{2(k+1)}}.
\end{aligned}
\end{equation}
Noticing that when $m>2-\frac{2}{N}$, it is easy to verify$$\frac{Q\lambda (k+m-1)}{2(k+1)}=\frac{(k+1)(m-2)+\lambda ^{*}(k+1)(3-m)}{(k+1)(k+m-1)\lambda ^{*}}< 1,$$
using Young's inequality, then
\begin{equation}\label{4.1.6}
\begin{aligned}
&C_{1}(N,k,m)(\left \| u^{\frac{k+m-1}{2}} \right \|_{L^{\frac{2(k+1)}{k+m-1}}(\mathbb{R}^{N})}^{\frac{2(k+1)}{k+m-1} }\left \| u^{\frac{k+m-1}{2}} \right \|_{L^{\frac{2}{k+m-1}}(\mathbb{R}^{N})}^{\frac{2}{k+m-1}})^{\frac{Q\lambda (k+m-1)}{2(k+1)}}\\
 &\leq k\left \| u^{\frac{k+m-1}{2}} \right \|_{L^{\frac{2(k+1)}{k+m-1}}(\mathbb{R}^{N})}^{\frac{2(k+1)}{k+m-1} }\left \| u^{\frac{k+m-1}{2}} \right \|_{L^{\frac{2}{k+m-1}}(\mathbb{R}^{N})}^{\frac{2}{k+m-1}}+C_{2}(N,k,m).
\end{aligned}
\end{equation}
Substitute \eqref{4.1.3}-\eqref{4.1.6} into \eqref{4.1.1} to get
\begin{equation}\label{4.1.7}
  _{0}^{C}\textrm{D}_{t}^{\alpha }\int _{\mathbb{R}^{N}}u^{k}dx+k\int _{\mathbb{R}^{N}}u^{k}dx\leq C_{2}(N,k,m),
\end{equation}
from Corollary \ref{corollary1}, we obtain
\begin{equation}\label{4.1.8}
  \int _{\mathbb{R}^{N}}u^{k}dx\leq \left \| u_{0} \right \|_{L^{k}(\mathbb{R}^{N})}^{k}+\frac{C_{2}(N,k,m)T^{\alpha }}{\alpha \Gamma },t\in [0,T].
\end{equation}
\noindent \textbf{Step 2. The $L^{\infty}$ estimates}

On account of the above arguments, our last task is to give the uniform boundedness of solution for any $t>0$.
Denote $q_{k}=2^{k}+2$, by taking $k=q_{k}$ in \eqref{4.1.1}, we have
\begin{equation}\label{4.2.1}
  \begin{aligned}
  &_{0}^{C}\textrm{D}_{t}^{\alpha }\int _{\mathbb{R}^{N}}u^{q_{k}}dx+\frac{4q_{k}(q_{k}-1)}{(m+q_{k}-1)^{2}}\left \| \bigtriangledown u^{\frac{m+q_{k}-1}{2}} \right \|_{L^{2}(\mathbb{R}^{N})}^{2}+q_{k}\int _{\mathbb{R}^{N}}u^{q_{k}+1}dx\int _{\mathbb{R}^{N}}udx\\
  &+q_{k}\int _{\mathbb{R}^{N}}u^{q_{k}}dx\leq q_{k}\int _{\mathbb{R}^{N}}u^{q_{k}+1}dx
  \end{aligned}
\end{equation}
armed with Lemma \ref{lemma333}, letting
$$v=\frac{m+q_ {k}-1}{2},q=\frac{2(q_{k}+1)}{m+q_{k}-1},r=\frac{2q_{k-1}}{m+q_{k}-1},c_{0}=c_{1}=\frac{1}{2q_{k}},$$
one has that for $N \geq 1$,
\begin{equation}\label{4.2.2}
\begin{aligned}
  \left \| u \right \|_{L^{q_{k}+1}(\mathbb{R}^{N})}^{q_{k}+1}&\leq C(N)c_{0}^{\frac{1}{\delta _{1}-1}}(\int _{\mathbb{R}^{N}}u^{q_{k}}dx)^{\gamma _{1}}
+\frac{1}{2q^{k}}\left \| \bigtriangledown  u^{\frac{m+q_{k}-1}{2}} \right \|_{L^{2}(\mathbb{R}^{N})}^{2}\\
&\qquad+\frac{1}{2q_{k}}\left \| u \right \|_{L^{m+q_{k}-1}(\mathbb{R}^{N})}^{m+q_{k}-1},
\end{aligned}
\end{equation}
where
$$\gamma _{1}=1+\frac{q_{k}+q_{k-1}+1}{q_{k-1}+\frac{p(m-2)}{p-2}}\leq 2,$$
$$\delta _{1}=\frac{(m+q_{k}-1)-2\frac{q_{k-1}}{p^{*}}}{q_{k}-q_{k-1}+1}=O(1).$$
Substituting \eqref{4.2.2} into \eqref{4.2.1} and with notice that $\frac{4q_{k}(q_{k}-1)}{(m+q_{k}-1)^{2}}\geq 2$. It follows
\begin{equation}\label{4.2.3}
  \begin{aligned}
  &_{0}^{C}\textrm{D}_{t}^{\alpha }\int _{\mathbb{R}^{N}}u^{q_{k}}dx+\frac{3}{2}\left \| \bigtriangledown u^{\frac{m+q_{k}-1}{2}} \right \|_{L^{2}(\mathbb{R}^{N})}^{2}+q_{k}\int _{\mathbb{R}^{N}}u^{q_{k}}dx\\
 &\qquad +q_{k}\int _{\mathbb{R}^{N}}u^{q_{k}+1}dx\int _{\mathbb{R}^{N}}udx\\
&\leq q_{k}\int _{\mathbb{R}^{N}}u^{q_{k}+1}dx\leq c_{1}(N)q_{k}^{\frac{\delta _{1}}{\delta _{1}-1}}(\int _{\mathbb{R}^{N}}u^{q_{k-1}}dx)^{\gamma _{1}}+\frac{1}{2}\left \| u \right \|_{L^{m+q_{k}-1}(\mathbb{R}^{N})}^{m+q_{k}-1}.
  \end{aligned}
\end{equation}
Applying Lemma \ref{lemma333} with $$v=u^ {\frac{m+q_{k}-1}{2}},q=2,\gamma=\frac{2q_{k-1}}{m+q_{k}-1},c_{0}=c_{1}=\frac{1}{2}$$ noticing $q_{k-1}=\frac{(q_{k}+1)+1}{2}$, and using Young's inequality, we obtain
\begin{equation}\label{4.2.4}
  \begin{aligned}
  &\frac{1}{2}\left \| u \right \|_{L^{m+q_{k}-1}(\mathbb{R}^{N})}^{m+q_{k}-1}=\frac{1}{2}\int _{\mathbb{R}^{N}}u^{m+q_{k}-1}dx\\
&\leq c_{2}(N)(\int _{\mathbb{R}^{N}}u^{q_{k-1}}dx)^{\gamma _{2}}+\frac{1}{2}\left \| \bigtriangledown u^{\frac{m+q_{k}-1}{2}} \right \|_{L^{2}(\mathbb{R}^{N})}^{2}\\
&\leq \frac{2}{3}\int _{\mathbb{R}^{N}}udx\int _{\mathbb{R}^{N}}u^{q_{k}+1}dx+c_{3}(N)+\frac{1}{2}\left \| \bigtriangledown u^{\frac{m+q_{k}-1}{2}} \right \|_{L^{2}(\mathbb{R}^{N})}^{2},
  \end{aligned}
\end{equation}
where $$\gamma _{2}=1+\frac{m+q_{k}-q_{k-1}-1}{q_{k-1}}< 2.$$
By summing up \eqref{4.2.3} and \eqref{4.2.4}, with the fact that $\gamma _{1}\leq 2$ and $\gamma _{2}<2$, we have
\begin{equation}
\begin{aligned}
&_{0}^{C}\textrm{D}_{t}^{\alpha }\int _{\mathbb{R}^{N}}u^{q_{k}}dx+q_{k}\int _{\mathbb{R}^{N}}u^{q_{k}}dx\leq c_{1}(N)q_{k}^{\frac{\delta _{1}}{\delta _{1}-1}}(\int _{\mathbb{R}^{N}}u^{q_{k-1}}dx)^{\gamma _{1}}+c_{3}(N)\\
&\leq max\left \{ c_{1} (N),c_{3}(N)\right \}q_{k}^{\frac{\delta _{1}}{\delta _{1}-1}}[(\int _{\mathbb{R}^{N}}u^{q_{k-1}}dx)^{\gamma 1}+1]\\
&\leq 2max\left \{ c_{1} (N),c_{3}(N)\right \}q_{k}^{\frac{\delta _{1}}{\delta _{1}-1}}max\left \{ (\int _{\mathbb{R}^{N}}u^{q_{k-1}}dx)^{2},1 \right \}.
\end{aligned}
\end{equation}
 Let $K_{0}=max\left \{ 1,\left \| u_{0} \right \|_{L^{1}(\mathbb{R}^{N})},\left \| u_{0} \right \| _{L^{\infty }(\mathbb{R}^{N})}\right \}$, we have the following inequality for initial data
\begin{equation}\label{4.2.5}
  \int _{\mathbb{R}^{N}}u_{0}^{q_{k}}dx\leq (max\left \{ \left \| u_{0} \right \|_{L^{1}(\mathbb{R}^{N})},\left \| u_{0} \right \| _{L^{\infty }(\mathbb{R}^{N})}\right \})^{q_{k}}\leq K_{0}^{q_{k}}.
\end{equation}
 Let $d_{0}=\frac{\delta _{1}}{\delta _{1}-1}$, it is easy to that $q_{k}^{d_{0}}=(2^{k}+2)^{d_{0}}\leq (2^{k}+2^{k+1})^{d_{0}}$. By taking $\bar{a}=max\left \{ c_{1}(N),c_{3}(N) \right \}3^{d_{0}}$ in the Lemma \ref{lemma444}, we obtain
\begin{equation}\label{4.2.6}
  \int u^{q_{k}}dx\leq (2\bar{a})^{2^{k}-1}2^{d_{0}(2^{k+1}-k-2)}max\left \{ \underset{t\geq 0}{sup}(\int _{\mathbb{R}^{N}}u^{q}dx)^{2^{k}},k_{0}^{q_{k}} \right \}\frac{T^{\alpha }}{\alpha \Gamma (\alpha )}.
\end{equation}
 Since $q_{k}=2^{k}+2$ and taking the power $\frac{1}{q_{k}}$ to both sides of \eqref{4.2.6}, then the boundedness of the solution $u(x,t)$ is obtained by passing to the limit $k\rightarrow\infty$
 \begin{equation}\label{4.2.7}
   \left \| u(x,t) \right \|_{L^{\infty }(\mathbb{R}^{N})}\leq 2\bar{a}2^{2d_{0}}max\left \{ sup_{t\geq 0}\int _{\mathbb{R}^{N}} u^{q_{0}}dx,K_{0}\right \}\frac{T^{\alpha }}{\alpha \Gamma (\alpha )}.
 \end{equation}
 On the other hand, by \eqref{4.1.8} with $q_{0}>2$, we know
 $$\int _{\mathbb{R}^{N}}u^{q_{0}}dx\leq \int _{\mathbb{R}^{N}}u^{3}dx\leq \left \| u_{0} \right \|_{L^{3}(\mathbb{R}^{N})}^{3}+\frac{cT^{\alpha }}{\alpha \Gamma (\alpha )}\leq K_{0}^{3}+\frac{cT^{\alpha }}{\alpha \Gamma (\alpha )}.$$
 Therefore we finally have
 \begin{equation}
   \left \| u(x,t) \right \|_{L^{\infty }(\mathbb{R}^{N})}\leq c(N,\left \| u_{0} \right \|_{L^{1}(\mathbb{R}^{N})},\left \| u_{0} \right \|_{L^{\infty }(\mathbb{R}^{N})},T^{\alpha })=M.
 \end{equation}
 \begin{remark}
   The method to prove the $L^{r}$ estimate to the $L^{\infty}$ estimate of the equation is also mentioned in \cite{2014Ultra}. This section mainly uses fractional differential inequalities, Gagliardao-Nirenberg inequalities and use the $\int _{\mathbb{R}^{N}}\left | \bigtriangledown u \right |^{2}dx$ and $\int _{\mathbb{R}^{N}}u^{k+1}dx\int _{\mathbb{R}^{N}}udx$ to control nonlinear terms $\int _{\mathbb{R}^{N}}u^{k+1}dx$. Therefore, the $L^{r}$ estimate is obtained, and if $k=q_{k}$ is estimated on $L^{\infty}$, the global boundedness of the solution for the NTFRDE is proved in any dimensional space.
 \end{remark}





\bibliographystyle{elsarticle-num}
\bibliography{Ref}
\end{document}